\documentclass[review]{elsarticle}

\usepackage[utf8]{inputenc}
\usepackage{amsmath}
\usepackage{amsfonts}
\usepackage{amssymb}
\usepackage{graphicx}
\usepackage{fullpage}
\usepackage[pdftex]{hyperref}
\usepackage{color}
\usepackage[titletoc,title]{appendix}
\usepackage[mathscr]{eucal}
\newtheorem{thm}{Theorem}[section]
\newtheorem{lem}[thm]{Lemma}
\newtheorem{prop}[thm]{Proposition}
\newtheorem{corol}[thm]{Corollary}
\newtheorem{defn}{Definition}[section]
\newtheorem{remark}{Remark}[section]
\newtheorem{assum}{Assumption}[section]
\newcommand{\proof}{\noindent {\bf Proof.}~}

\DeclareMathOperator*{\diam}{diam}

\DeclareMathOperator*{\diag}{diag}
\DeclareMathOperator*{\interior}{int}

\newcommand{\bN}{\mathbb{N}}
\newcommand{\bR}{\mathbb{R}}
\newcommand{\bZ}{\mathbb{Z}}

\newcommand{\bu}{\mathbf{u}}
\newcommand{\bv}{\mathbf{v}}
\newcommand{\bx}{\mathbf{x}}

\newcommand{\bA}{\mathbf{A}}

\newcommand{\cI}{\mathcal{I}}
\newcommand{\cN}{\mathcal{N}}
\newcommand{\cU}{\mathcal{U}}
\newcommand{\cV}{\mathcal{V}}

\newcommand{\powerset}{\raisebox{.15\baselineskip}{\Large\ensuremath{\wp}}}

\begin{document}

\title{The Echo Index and multistability in input-driven recurrent neural networks}

\author[1]{Andrea Ceni\corref{cor1}}\ead{ac860@exeter.ac.uk}

\author[2]{Peter Ashwin}\ead{p.ashwin@exeter.ac.uk}

\author[3,1]{Lorenzo Livi}\ead{lorenzo.livi@umanitoba.ca}

\author[4]{Claire Postlethwaite}\ead{c.postlethwaite@auckland.ac.nz}

\cortext[cor1]{Corresponding author}
\address[1]{Department of Computer Science, University of Exeter, Exeter EX4 4QF, UK}
\address[2]{Department of Mathematics, University of Exeter, Exeter EX4 4QF, UK}
\address[3]{Departments of Computer Science and Mathematics, University of Manitoba, Winnipeg, MB R3T 2N2, Canada}
\address[4]{Department of Mathematics, University of Auckland, Auckland 1142, New Zealand}

\begin{abstract}
A recurrent neural network (RNN) possesses the echo state property (ESP) if, for a given input sequence, it ``forgets'' any internal states of the driven (nonautonomous) system and asymptotically follows a unique, possibly complex trajectory. The lack of ESP is conventionally understood as a lack of reliable behaviour in RNNs. Here, we show that RNNs can reliably perform computations under a more general principle that accounts only for their local behaviour in phase space.
To this end, we formulate a generalisation of the ESP and introduce an echo index to characterise the number of simultaneously stable responses of a driven RNN.
We show that it is possible for the echo index to change with inputs, highlighting a potential source of computational errors in RNNs due to characteristics of the inputs driving the dynamics.
\end{abstract}
\begin{keyword}
Nonautonomous dynamical systems \sep Input-driven systems \sep Recurrent neural networks \sep Echo state property \sep Multistability \sep Machine learning.
\end{keyword}

\maketitle

\newpage

\tableofcontents

\newpage

\section{Introduction}

Recurrent neural networks (RNNs) are input-driven (i.e. nonautonomous) dynamical systems \cite{manjunath2012theory} whose behaviour depends both on model parameters and on inputs to the system. In order to describe responses of a (trained) RNN to the full range of possible inputs, it is necessary to go beyond the theory of autonomous (input-free) dynamical systems and consider more general nonautonomous dynamical systems \cite{kloeden2011nonautonomous}, where the equations ruling the dynamics change over time. The theory of nonautonomous dynamical systems is much less developed than that of autonomous systems, and notions such as convergence (and hence attractors) need to be carefully defined \cite{kloeden2016forward}.

Starting with the work of Jaeger \cite{jaeger2001echo}, several authors have proposed that a successfully trained RNN should have the so-called \emph{echo state property} (ESP) \cite{yildiz2012re}. The idea of ``echo state'' gave rise to a training paradigm of RNNs called \emph{reservoir computing} \cite{lukovsevivcius2009reservoir} and a class of RNNs known as \emph{echo state networks} (ESNs) \cite{jaeger2004harnessing}.
An ESN is an RNN that is relatively easy to train, since optimisation is restricted to the output layer and recurrent connections are left untouched after initialisation.

If an RNN possesses the ESP this means that, for a given input sequence, it will asymptotically produce the same sequence of states and will ``forget'' any internal state, ending up following a unique (though possibly very complex) trajectory in response to that input \cite{manjunath2013echo}.
This trajectory represents the solution of the specific problem encoded in the input sequence that is fed to the network and so the system can be seen as acting as a \emph{filter} that transforms the input sequence into a unique sequence of output \cite{grigoryeva2018echo}.
The presence of the ESP has been historically associated with reliable RNN behaviour. Accordingly, the loss of ESP has been directly associated with the loss of reliable behaviour in RNNs, implying that correct computation is not possible without ESP.

In this paper, we analyse the ESP through the lens of nonautonomous dynamical system theory and introduce a generalisation of the ESP that accounts only for local behaviour in phase space.
To this end, we introduce an ``echo index'' to characterise the number of simultaneously stable responses of a driven RNN. 
We achieve this by defining a simple, yet non-trivial class of attracting solutions for input-driven systems that we call \emph{uniformly attracting entire solutions} and count how many of these solutions coexist in phase space under the action of an input sequence.
We show that echo index 1 and the classical notion of ESP do not always coincide.
More precisely, driven RNNs having echo index 1 are allowed to produce reliable behaviour, in the ESP spirit, only in a phase space subset thus eliminating those regions producing unreliable behaviour but having zero probability of occurring.

The presence of multiple attractive solutions may be useful when the RNN is used to perform context-dependent computations. For example, consider the addition modulo $ N $, i.e. the addition in a modular arithmetic. Given an initial integer value $ n_0 $ and an input sequence of integers $ n_1 , n_2, \ldots $, assuming values in $\{0,1, \dots, N-1\}$, the neural network needs to compute online the sum $ N_k:= \sum_{i=0}^{k} n_i \,\, (\text{mod } N )$ as time $k$ runs. Therefore, in this case the system has to respond in $ N $ different ways depending on the initial state.

The remainder of this paper is structured as follows. 
Background material is introduced in Section \ref{sec:nds}. 
In Section \ref{sec:nesp}, we describe our main contribution: a generalisation of the ESP for driven RNNs.
Our main theoretical results are discussed in Section \ref{sec:theoretical_results}: (1) in Theorem \ref{thm:idea_contraction} we provide a sufficient condition for the existence and uniqueness of a uniformly attracting entire solution in a phase space subset; (2) in Theorem \ref{thm:small_input}, we discuss a mechanism giving rise to the occurrence of multistability in systems driven by low amplitude inputs; (3) in Proposition \ref{prop:large_inputs} we prove that forcing RNNs with large-amplitude inputs will induce echo index 1 (and the ESP) and finally, (4) in Section~\ref{sec:input-stable} we discuss stability of the echo index w.r.t. perturbations on the input sequences and show, in Theorem \ref{thm:dense}, how this depends on the metric imposed on the space of input sequences.
In Section \ref{sec:experiments} we report some numerical experiments.
Section \ref{sec:2esp} provides an example of multistable RNN dynamics possessing echo index 2.
In Section \ref{sec:Change-echo-index}, we show how the echo index might depend, fixing the RNN model parameters, on the particular input driving the dynamics, illustrating a bifurcation from echo index 1 to echo index 2.
This offers new insights on the fact that changes in behaviour may also occur for reasons that are not related to changes in model parameters (i.e. via training).
In Section \ref{sec:Maass}, we train an RNN to solve the context-dependent task illustrated in \cite[Figure 5]{hoerzer2012emergence} and give an interpretation of the related multistable, nonautonomous RNN dynamics based on the modeling framework developed here.
Finally, Section \ref{sec:conclusions} concludes the paper and provides final remarks.
Some background definitions, results and proofs are provided in appendices located at the end of the manuscript. In particular, we show in Appendix \ref{sec:pullback_thms} (Proposition \ref{prop:pullback_thm}) that the \emph{natural association} of \cite{manjunath2013echo} corresponds to the system's pullback attractor \cite{kloeden2011nonautonomous} and in Appendix \ref{sec:forward} we link pullback convergence with uniform forward convergence in RNNs.

\section{Nonautonomous dynamics of recurrent neural networks}
\label{sec:nds}

In this section, we investigate RNN dynamics using tools from nonautonomous dynamical systems theory.
In Section~\ref{sec:rec_neur_net}, we introduce RNNs and provide details about the RNN model considered in our analysis.
Then, in Section~\ref{sec:inp-driv-DS}, we highlight the input-driven nature of RNN dynamics and, successively in Section~\ref{sec:cocycle_formalism}, we introduce the skew product formalism for nonautonomous dynamical systems \cite{kloeden2011nonautonomous}.
Finally, in Section~\ref{sec:pullback_attractors}, we introduce the fundamental notion of \emph{pullback attractor}.

\subsection{Recurrent neural networks}
\label{sec:rec_neur_net}

A discrete-time RNN \cite{bianchi2017recurrent} is a state-space model describing the evolution of state variable $x[k]\in X\subset \bR^{N_r}$ in discrete time $k\in\bZ$.
In this paper, we consider a fairly general RNN model equipped with leaky-integrator neurons \cite{jaeger2007optimization},
\begin{align}
\label{eq:leaky_rnn}
x[k+1] &  = (1-\alpha) x[k] + \alpha \phi( W_r x[k] + W_{in} u[k+1] + W_{fb} z[k] ), \\
\label{eq:output}
z[k+1] & = \psi(x[k+1]).
\end{align}
$\phi(\cdot)$ is a component-wise activation function (e.g. hyperbolic tangent) and $u[k]\in U$ is the $N_i$-dimensional input sequence for $U\subset \bR^{N_i}$. The sets $X$ and $U$ denote the state (or phase) and input spaces, respectively. The matrices $W_r \in \mathbb{R}^{N_r \times N_r} $ and $ W_{in} \in \mathbb{R}^{N_r \times N_i} $ represent recurrent and input-to-network couplings. The output feedback matrix $W_{fb}\in\mathbb{R}^{N_r\times N_o}$ injects the last computed output into the state-update equation \eqref{eq:leaky_rnn}.
The scalar $\alpha \in (0,1] $ can be used as hyperparameter to control the RNN time-scale given by $1/\alpha$ \cite{tallec2018can,jaeger2007optimization}.
The particular implementation of $\psi(\cdot)$ in \eqref{eq:output} depends on the task at hand. For instance, in classification tasks $\psi(\cdot)$ might take the form of a softmax assigning probabilities to predicted classes; in forecasting the easiest choice is a linear deterministic function, i.e. $z[k+1] = \psi(x[k+1]) = W_o x[k+1]$, although nonlinear functions are also common.

Training an RNN means to optimise its parameters, e.g. the entries of $ W_r, W_{in}, W_{fb}, W_o $ and $\alpha$, in order to achieve a suitable configuration of the system for solving a given task.
This is typically accomplished by means of a gradient descent algorithm or variation of thereof \cite{ruder2016overview}.
Learning long-term dependencies with gradient descent is known to be problematic, as a consequence of the so-called vanishing/exploding gradient problem \cite{pascanu2013difficulty}. For this purpose, two types of methods have been proposed: (i) gating mechanisms \cite{hochreiter1997long,chung2014empirical} and (ii) approaches based on unitary matrices and constant-slope activation functions \cite{NIPS2016_6327}.
ESNs \cite{7386673,jaeger2004harnessing,pascanu2011neurodynamical,tivno2018asymptotic,esnfish2016,grigoryeva2018echo}, a special class of RNNs, bypass this problem as training targets the output layer weights $W_o$ only. The recurrent layer, called a reservoir, is randomly instantiated (although more sophisticated methods have been proposed in the literature \cite{rodan2012simple}) and the model is modified only offline at the hyper-parameter level.
This simple training protocol is not sufficient in many applications, e.g. when it is required to learn memory states. To this end, training mechanisms based on output feedback \cite{8720322,PhysRevLett.118.258101} and online training \cite{hoerzer2012emergence,sussillo2009generating} have been proposed, with successful applications in physics \cite{seoane2019evolutionary,lu2018attractor}, complex systems modeling \cite{ibanez2018detection,PhysRevE.98.052209}, and neuroscience \cite{buonomano2009state}, just to name a few.

In the remainder of this paper, we will assume to deal with an already trained RNN and thus do not consider the effects of training on the dynamics.

\subsection{Input-driven dynamical systems}
\label{sec:inp-driv-DS}

The equations \eqref{eq:leaky_rnn}-\eqref{eq:output} ruling the behaviour of a (trained) RNN can be seen as a special case of an input-driven dynamical system of the form:
\begin{equation}
\label{eq:rnn_nds}
x[k+1] = G(u[k+1], x[k]),
\end{equation}
where the map $ G:U \times X \rightarrow X$ is defined as
\begin{equation}
\label{eq:RNN-map}
    G(u,x) = (1-\alpha)x + \alpha \phi( W_r x + W_{in} u + W_{fb} \psi(x) ).
\end{equation}
The action of inputs $u[k]$ driving the dynamics of $x[k]$ gives an explicit time-dependence of \eqref{eq:leaky_rnn} and means this system is nonautonomous. A given input sequence $ \bu = \{ u[k] \}_{k \in \bZ}$ induces a sequence of maps $ \{ f_k \}_{k \in \bZ} $, where $ f_k(\cdot) := G(u[k+1], \cdot ): X \rightarrow X $ is the map ruling the update of the RNN state at time $k$, i.e. $x[k+1]=f_k(x[k])$.

The first-order derivative matrix of $G$ \eqref{eq:RNN-map} w.r.t. the state variable reads
\begin{equation}
\label{eq:jacobian_G}
D_x G(u,x) = (1-\alpha)I_{N_r} +  \alpha S(u,x) M(x),
\end{equation}
where
\begin{equation}
\label{eq:effective_recurrent_matrix}
M(x)= W_r + W_{fb}  D_x \psi(x)
\end{equation}
denotes the ``effective recurrent matrix'', $ D_x \psi(\cdot) $ is the Jacobian matrix of $\psi$ \eqref{eq:output}, and
\begin{equation}
\label{eq:diag_deriv}
S(u,x)= \diag\left[ \phi' \left( \xi_j(u,x) \right) \right]_{j=1,\ldots,N_r},
\end{equation}
with pre-activations $\xi_j(u,x):= (W_r)_{(j)}\cdot x + (W_{fb})_{(j)} \cdot \psi(x) + (W_{in})_{(j)} \cdot u$.
Note that, when considering the map \eqref{eq:RNN-map}, as long as the neuronal activation function $ \phi $ and the readout function $ \psi $ are both regular of class $C^1$, then the map $G$ will be too.
Moreover, for all $u \in U$, the map $G(u,\cdot): X \rightarrow X $ is a local diffeomorphism whenever matrix $M(x)$ \eqref{eq:effective_recurrent_matrix} is invertible.
In addition, if $\phi$ is bounded with image $(-L,L)$ then the state space $X$ can be assumed to be the compact hypercube $[-L,L]^{N_r} := \{ (x_1, \ldots, x_{N_r} ) \in \bR^{N_r} \,\,|\,\, x_i \in [-L,L], \, i= 1, \ldots, N_r  \} $, e.g. if $ \phi = \tanh $ then $L=1$, see Proposition \ref{prop:absorbing-set-pullback}.

We make the following standing assumptions:
\begin{assum}
\label{assum:fundamental}
\begin{itemize}
    \item[(i)] $G$ is continuously differentiable in all arguments, i.e $G \in C^1( U\times X , X )$;
    \item[(ii)] for all $ u \in U $, the map $ G(u,\cdot): X \rightarrow X $ is a local diffeomorphism onto its image;
    \item[(iii)] $ U \subset \bR^{N_i} $ is compact and $ X \subset \bR^{N_r} $ is usually the compact closure of a $N_r$-dimensional Cartesian product of real intervals.
\end{itemize}
\end{assum}
\begin{remark}
(ii) implies that the preimage of any zero measure set is also zero measure, hence given any input sequence $ \{ u[k] \}_{k \in \bZ} $ assuming values in $U$, for any $Z\subset X$ with $\lambda(Z)=0$, then $\lambda(f_k^{-1}(Z))=0$ for all $k \in \bZ$, where $\lambda $ denotes Lebesgue measure on $X\subset \bR^{N_r}$. This is weaker than assuming that $G(u,\cdot)$ is invertible for fixed $u$, but it means that phase space volume cannot ``suddenly collapse''.
\end{remark}

\subsection{The cocycle formalism}
\label{sec:cocycle_formalism}

There are two ways to describe nonautonomous systems \cite{kloeden2011nonautonomous}: the \emph{process} and the \emph{skew product} (also called \emph{cocycle}) formalism. In this paper, we use the cocycle formalism that is convenient when describing the input evolution as a shift in sequence space.

Let $(X, d_X) $ and $(U,d_U)$ be compact metric spaces.
Time evolution will be parametrised through the ring of integers $\bZ$ or a subset of it. We write $ \bZ^+:=\{k\in\bZ~:~k\geq 1\}$ and $\bZ_{0}^-:=\{k\in\bZ~:~k\leq 0\}$. 
Let $\mathbb{T}$ be one of the sets $\bZ, \,\bZ^+,  $ or $\bZ_{0}^-$, we consider the set $U^{\mathbb{T}}:=\{ {\bf u}=\{u[k]\}_{k\in \mathbb{T}}~:~u[k]\in U, \,\, \forall k \in \mathbb{T}\}$ of all input sequences assuming values in set $U$, and we will denote $\cU=U^{\bZ}$, $\cU^+=U^{\bZ_{0}^+}$ and $\cU^-=U^{\bZ_{0}^-}$.
Moreover, given ${ \bf u } \in \cU$, we will denote with $ { \bf u }^+$ and ${ \bf u }^-$ the projection of $\bf u$ to $\cU^+,$ and $\cU^-$, respectively.
The set $\cU$ is usually equipped with the product topology induced by the metric
\begin{equation}
    \label{eq:metric-inputsequences}
    d_{\text{prod}}({\bf u},{\bf v}):= \sum_{k \in \mathbb{Z}}\dfrac{d_U(u[k],v[k])}{2^{|k|}},
\end{equation}
which renders $(\cU, d_\text{prod})$ a compact metric space.
Another metric suitable for our purposes is the uniform metric, defined as:
\begin{equation}
    \label{eq:uniform-metric-inputsequences}
    d_{\text{unif}}({\bf u},{\bf v}):= \sup_{k \in \mathbb{Z}}{d_U(u[k],v[k])}.
\end{equation}
From the fact that $U$ is compact, it follows that $(\cU, d_\text{unif})$ is a compact metric space.

The dynamics on this space of input sequences is described by means of the \emph{shift operator}, which is a map $\sigma: \cU \rightarrow \cU$ defined as follows:
\begin{equation}
\label{eq:shift}
    {\bf u}:=\{ u[k]\}_{k \in \mathbb{Z}} \longmapsto \sigma({\bf u}):= \{ u[k+1]\}_{k \in \mathbb{Z}}.
\end{equation}
The composition $\sigma^n$ defines a discrete dynamical system on the metric space $(\cU,d_{\cU})$ \cite{kloeden2011nonautonomous}. Setting $\sigma^0( {\bf u} ) := {\bf u}$ we have that $\{ \sigma^n \}_{n \in \mathbb{Z}}$ represents a group of homeomorphisms on $\cU$, which expresses the sequential forward or backward shift in time of all input sequences.
Moreover, defining $ p: \cU \rightarrow U $ as the projection mapping
\begin{equation}
\label{eq:projection}
    {\bf u}:=\{ u[k]\}_{k \in \mathbb{Z}} \longmapsto p({\bf u}):=  u[0],
\end{equation} 
we get in $U$ the current value of input sequence $ \bf u$ as $ p(\sigma^n({\bf u}))= u[n]$.

We describe the dynamics in response to input using a cocycle map \cite[Definition 2.1, page 28]{kloeden2011nonautonomous} for RNNs.

\begin{defn}
\label{def:cocycle_mapping}
The nonautonomous dynamical system \eqref{eq:rnn_nds} can be described using a \emph{cocycle mapping}, $\Phi~:~ \mathbb{Z}_0^+ ~\times~ \cU~ \times~ X \longrightarrow~ X $, defined as follows:
\begin{align}
\label{eq:cocycle1}
    \Phi(0, {\bf u}, x_0 ) &:= x_0, &\forall {\bf u} \in \cU, \forall x_0 \in  X, \\
    \label{eq:cocycle2}
    \Phi(n, {\bf u}, x_0 ) &:= G( p(\sigma^n({\bf u})),  \Phi(n-1, {\bf u}, x_0 )  ), & \forall n \geq 1,  \forall {\bf u} \in \cU, \forall x_0 \in  X
\end{align}
\end{defn}

Definition \ref{def:cocycle_mapping} leads to the following \eqref{eq:cocycle_property} to hold. We state the result without proof since it is well-known in the literature under the name of \emph{cocycle property}. 
\begin{lem}
\label{lem:cocycle_property}
The set $\{ \Phi(n,*,\cdot) \}_{n \in \mathbb{Z}_0^+}$ forms a semigroup of continuous functions from $\cU \times X $ to $ X $.
In particular, relations \eqref{eq:cocycle1}-\eqref{eq:cocycle2} imply the \emph{cocycle property}:
\begin{equation}
\label{eq:cocycle_property}
   \Phi(m+n, {\bf u}, x_0 )  \,\,  =  \,\,  \Phi(n, \sigma^m({\bf u}), \Phi(m, {\bf u}, x_0 ) ) ,
\end{equation}
for any $m,n\in\bZ_0^+$, $x_0 \in X$ and $\bu \in \cU$.
\end{lem}

Note that \eqref{eq:rnn_nds} implies that the forward map is always defined, i.e. given a point $x_0 \in X$ and any input ${\bf u} \in \cU$, the forward trajectory is uniquely defined by \eqref{eq:rnn_nds}. On the other hand, it is possible that a backward trajectory does not exist, or may not be unique if one does exist.
For example, although the one-dimensional map $ G(u,x) :=  \tanh( \mu x + u ) $ is invertible for any fixed $\mu \in (0,1)$, every backward trajectory constructed from a $x\neq 0$ leads outside the compact set $[-1,1]$ in a finite number of backward steps, thus making impossible to obtain a further preimage of $\tanh$.
This means it may not be possible to extend the cocycle mapping backward in time.
Trajectories that are well-defined in the infinite past play an important role in the nonautonomous dynamics~\cite{manjunath2013echo}, as expressed in the next definition.\footnote{We do not consider invariant sets or entire solutions in terms of the skew product formalism \cite[Definition 2.19]{kloeden2011nonautonomous}, but rather in terms of the process \cite[Definition 2.14]{kloeden2011nonautonomous} induced by a given input sequence $\bu$.}
\begin{defn}
\label{def:entire_solution}
An \emph{entire solution} for the system in \eqref{eq:rnn_nds} with input $\bu:=\{ u[k] \}_{k \in \bZ} \in \cU$ is a bi-infinite sequence of states $\{ x[k] \}_{k \in \bZ}$ that satisfies \eqref{eq:rnn_nds} for all $k \in \bZ$. 
In other words,
$$
\Phi(s, \sigma^m(\bu),  x[m]) =   x[m+s]
$$
for all $m \in \bZ$ and  $s \in \bZ_0^+$. 
\end{defn}

Assuming the existence of an entire solution $\{ x[k] \}_{k \in \bZ}$ for input $\bu \in \cU$, and exploiting the forward definition of the cocycle mapping, we can write the past evolution as follows:
\begin{equation}
\label{eq:cocycle_back}
    x[-n] = \Phi(m-n, \sigma^{-m}(\bu), x[-m]),  \quad \forall n \in \bZ_0^+ , \, \forall m\in\bZ_0^+  \mbox{ with }  m \geq n.
\end{equation}
Such a relation expresses the fact that the point $x[-n]$ is the resulting state of the system if we start from $x[-m]$ and drive the dynamics with the sequence of input values $u[-m+1],\ldots, u[-n]$.

\subsection{Pullback attractors}
\label{sec:pullback_attractors}

In this section we introduce some basic definitions of the theory of nonautonomous dynamical systems, including the one of pullback attractor. In fact, the ESP notion first introduced in \cite{jaeger2001echo} was formulated using left-infinite sequences and a pullback argument.
The following definition extends the notion of entire solution.
\begin{defn}
\label{def:invariant_nonaut_family}
Consider a nonautonomous system defined by a cocycle mapping as in Definition \ref{def:cocycle_mapping} with input sequence $\bu\in\cU$.
A family of nonempty compact sets $\bA=\{ A_n \}_{n\in \bZ}$ is called an \emph{invariant nonautonomous set} for input $\bu$ if 
$$
\Phi(s, \sigma^{m}(\bu), A_m) =  A_{s+m} .
$$
for all $m \in \bZ$ and $s \in \bZ_0^+$.
\end{defn}

Entire solutions are invariant nonautonomous sets where each $A_n$ is a single point. Invariant nonautonomous sets turn out to be composed by entire solutions \cite[Lemma 2.15]{kloeden2011nonautonomous}.
Replacing ``$ = $'' with ``$ \subseteq $'' in Definition \ref{def:invariant_nonaut_family}, we obtain the definition of a positively invariant family of sets.
\begin{defn}
\label{def:posit_invar_family}
A family of nonempty compact sets $\mathbf{B}=\{ B_n \}_{n\in \bZ}$ is called a \emph{positively invariant nonautonomous set} for input $\bu$ (or simply $\bu$\emph{-positively invariant}) if 
$$
\Phi(s, \sigma^{m}(\bu), B_m) \subseteq  B_{s+m} .
$$
for all $m \in \bZ$ and $s \in \bZ_0^+$.
\end{defn}

Note that the assumption that $G$ is well-defined as a map to $X$ implies that $X$ itself is positively invariant.
Invariant sets play an important role for understanding the behaviour of a dynamical system. 
The most relevant ones are those invariant sets that attract the surrounding trajectories.
Since, in nonautonomous systems, the equations ruling the dynamics change with time, the notion of attraction can be formulated in various ways, leading to definitions of \emph{forward attraction} and \emph{pullback attraction}. It is interesting to note that pullback attractors share more properties with their autonomous counterparts  \cite{kloeden2016forward}.

Based on \cite[Definitions 3.3 and 3.4]{kloeden2011nonautonomous}, we introduce the notion of (global) pullback attractor as follows.\footnote{Note that as $X$ is bounded, all bounded subsets of $X$ are pullback attracted to $A_n$ \cite[Definition 3.4]{kloeden2011nonautonomous} if and only if $X$ is.}
Let $h$ be Hausdorff semi-distance for the metric space $(X,d_{X})$; see Appendix \ref{sec:hausdorff} for details.
\begin{defn} 
\label{def:pullback_attractor}
An invariant nonautonomous set $\bA=\{ A_n \}_{n\in \bZ}$ for input $\bu$ consisting of nonempty compact sets is called a \emph{(global) pullback attractor} of the system \eqref{eq:rnn_nds} driven by $\bu$ if
\begin{equation}
    \label{eq:pullback_attractor}
        \lim_{k \rightarrow \infty} h( \Phi(k, \sigma^{-k+n} (\bu), X ), A_n ) = 0, \qquad \forall n \in \bZ.
\end{equation}
\end{defn}

It is easy to show from this that if $\bA=\{A_n\}_{n \in \bZ}$ is a pullback attractor for input $\bu$ then, for any $m\in\bZ$, $\{A_{n+m}\}_{n \in \bZ}$ is a pullback attractor for $\sigma^m(\bu)$.

Manjunath and Jaeger~\cite{manjunath2013echo} introduce the notion of \emph{natural association} as follows.
\begin{defn} \cite[Definition 4]{manjunath2013echo}
\label{def:natural_association}
Consider an input sequence $\bv:=\{ v[k] \}_{k \in \bZ} \in \cU$.
The sequence of sets $\{ X_n(\bv) \}_{n \in \bZ}$ defined by 
\begin{equation}
    \label{eq:natural_association}
    X_n(\bv) : = \bigcap_{m \leq n} \Phi( n-m , \sigma^{m}(\bv), X), \qquad  n \in \bZ,
\end{equation}
is called the \emph{natural association} of the process induced by input sequence $\bv$.
\end{defn}
\begin{remark}
In Proposition \ref{prop:pullback_thm}, we rigorously prove that the nonautonomous set $ \{ X_n ( \mathbf{v} ) \}_{n \in \bZ}  $ of Definition \ref{def:natural_association} is the unique invariant compact nonautonomous set $ \{ A_n \}_{n \in \bZ} $ which satisfies equation \eqref{eq:pullback_attractor} for the input sequence $\bv$.
For this reason, from now on we will call the nonautonomous set $ \{ X_n ( \mathbf{u} ) \}_{n \in \bZ}  $ of \eqref{eq:natural_association} the global pullback attractor for input $\bu$.
\end{remark}

\section{An echo index for recurrent neural networks} 
\label{sec:nesp}

We first consider the ESP for driven RNNs as defined in \cite{manjunath2013echo}. Then, we introduce a more general definition accounting for reliable RNN behaviour in a phase space subset only. This generalisation also leads to a well-defined notion of multistable behaviour for driven RNNs.
\subsection{Entire solutions and the Echo State Property}
\label{sec:entire_pullback}

Recall from Definition \ref{def:entire_solution} that an {\em entire solution} for input $\mathbf{v}=\{v[k]\}_{k\in\bZ}\in \cU$ is a sequence $\mathbf{x}=\{x[k]\}_{k\in\bZ}$ that satisfies \eqref{eq:rnn_nds}, for all $k\in\bZ$.
\begin{defn} \cite[Definition 2]{manjunath2013echo}
\label{def:uniqueentiresolution-ESP}
The nonautonomous system \eqref{eq:rnn_nds} has the \emph{Echo State Property} for input $\bv \in \cU$ if there exists a unique entire solution $\bx$. 
\end{defn}

We give a simple proof of the fact that, if the ESP holds for a given $\bv \in \cU$, then the entire solution is the pullback attractor as in Definition \ref{def:pullback_attractor}, using the notation of the cocycle mapping introduced in Section \ref{sec:cocycle_formalism}.
\begin{prop}
\label{prop:ESP_pullback}
If \eqref{eq:rnn_nds} has the ESP for input $\bv \in \cU$, then the (unique) entire solution is the global pullback attractor of the system \eqref{eq:rnn_nds} with input $\bv$.
\end{prop}

\proof
First, note that if $\{\bx^*[k]\}_{k \in \bZ}$ is an entire solution for input $\bv \in \cU$, then by the definition of $\Phi$, it holds true that $ x^*[n] = \Phi(n-m, \sigma^{m}(\bv), x^*[m])\,\, \forall \,m < n , \, \forall n \in \bZ$. Compactness of $X$ guarantees that $ x^*[m] \in X$ from which it follows that $ x^*[n] \in \bigcap_{m<n} \Phi(n-m, \sigma^m(\bv), X) \,\, \forall n \in \bZ $.
The ESP hypothesis means $\{\bx^*[k]\}_{k \in \bZ}$ is the only entire solution for input $\bv \in \cU$. 
Therefore, we need to prove that for all $n \in \bZ$ the pullback attractor is composed of the singleton $ x^*[n] = \bigcap_{m<n} \Phi(n-m, \sigma^{m}(\bv), X)$.
We prove this by contradiction: suppose there exists a $n_0 \in \bZ$, where w.l.o.g. $n_0=0$, and a $ y_0 \in X$ with $y_0 \neq x^*[0]$ such that
$$
\{ x^*[0], y_0 \} \subset \bigcap_{m<0} \Phi(-m, \sigma^{m}(\bv), X).
$$
This means that $y_0 \in \Phi(-m, \sigma^{m}(\bv), X) \,\, \forall m\leq -1 $, i.e. for each $m\leq -1$ there exists at least one point, let us denote it $y[m] \in X$, such that $ y_0 = \Phi(-m, \sigma^{m}(\bv), y[m])$. In other words, there exists a backward trajectory $ \mathbf{y}^{-}:= \{ \ldots, y[-2], y[-1], y[0]:=y_0 \}$ for the point $y_0$. Now, by iterating \eqref{eq:rnn_nds} forward in time from $y[0]:=y_0$, we get the entire solution $ \mathbf{y}:=\{ y[k]\}_{k \in \bZ}$. This gives a contradiction.
\qed\\

Jaeger \cite{jaeger2001echo} shows that, if the ESP holds \emph{for all input sequences} assuming values in a given compact set, then the unique entire solution is forward attracting \cite[Definition 4]{jaeger2001echo}.
In Appendix \ref{sec:forward}, Proposition \ref{prop:ESP_unif-forw-attract}, we provide an alternative proof of this fact by applying \cite[Theorem 3.44]{kloeden2011nonautonomous} to RNNs as special class of nonautonomous dynamical systems.

\subsubsection{Limitations of pullback attractors in describing reliable RNN behaviour}
\label{sec:kloeden}

Here, we emphasise the difference between pullback and forward mechanisms of attraction with an example of Kloeden \textit{et al.} \cite{kloeden2012limitations}. 

\begin{figure}[ht!]
\begin{center}
\includegraphics[keepaspectratio=true,scale=0.31]{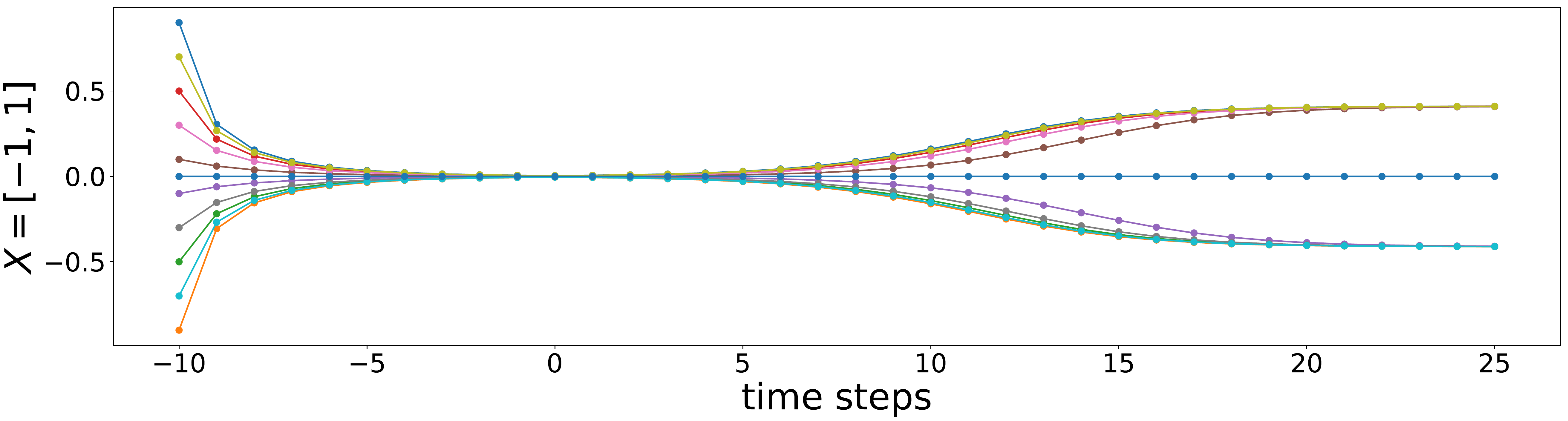}
\end{center}
\caption{
Equation \eqref{eq:kloeden}, with $a=1.5$, has been evolved from $k=-10$ to $k=25$ starting with 11 initial conditions in $[-1,1]$. One of these initial conditions is $x=0$ which is exactly on the global pullback attractor of the system driven by the input sequence $u[k] = a$ for $k\geq0$ and $u[k]=a^{-1}$ for $k<0$, with $a=1.5$. 
}
\label{fig:kloeden}
\end{figure}

Let us consider an input-driven system where $X$ is the compact $[-1,1]$ and the equation ruling the dynamics is
\begin{equation}
    \label{eq:kloeden}
        x[k+1] = \tanh\left( u[k] \dfrac{x[k]}{1+\left|x[k]\right|} \right) .
\end{equation}
Now, let us drive the system with the input sequence $u[k] = a$ when $k\geq0$, and $u[k]=a^{-1}$ when $k<0$, for some constant $a>1$.
The process generated by this input sequence has only one entire solution, i.e. $x[k]=0$, for all $k \in \mathbb{Z}$.
Accordingly, Proposition~\ref{prop:ESP_pullback} implies that such an entire solution $ x[k]\equiv 0 $ is the global pullback attractor.
Nevertheless, that unique entire solution is ``attractive'' in the past, yet it is ``repulsive'' in the future. Moreover, all initial conditions (except $x=0$) asymptotically tend to assume two possible values forward in time; see Figure \ref{fig:kloeden} for a visual representation.

Definition~\ref{def:uniqueentiresolution-ESP} implies that the system taken into account here has the ESP and therefore we would conclude that it produces a unique reliable response. However, our example shows how such a conclusion might be misleading.
Furthermore, the unique entire solution which is supposed to be the ``echo'' of the input sequence is actually produced by the system exclusively when the system is initialised exactly at $x_0=0$.

It appears clear that the mere consideration of pullback attraction w.r.t. an input sequence is not enough to meet the ideas behind the ESP (i.e. the idea of reliable, unique asymptotic behaviour in response to an input sequence).

\subsection{Uniformly attracting entire solutions and the echo index}
\label{sec:localPullback_echoIndex}

If the ESP (Definition~\ref{def:uniqueentiresolution-ESP}) does not hold, then clearly a wide variety of behaviours are possible.
Here, we are interested in describing the case where there is a finite number of stable responses to an input sequence.
For this we state a notion of \emph{uniformly attracting entire solution} (UAES).\footnote{Compared to \cite[Definition 3.48(iii)]{kloeden2011nonautonomous}: although the attraction is uniform, we do not require the neighbourhood to be uniform in $k$. On the other hand, this is a special case that only considers entire solutions that are attractors.} This is a local attractor: as noted previously \cite{Ochs1999}, it is possible for a pullback attractor to have a decomposition into a number of local attractors. 
\begin{defn}
\label{def:local_point_attractor}
Consider a fixed input sequence $\bu \in \cU$, an entire solution $ \{ x[k] \}_{k \in \bZ}$ and a positively invariant nonautonomous set $\{ B[k] \}_{k \in \bZ}$ composed of compact sets.
\begin{itemize}
    \item[(i)] If
    \begin{equation}\label{eq:uni_pt_att}
        \lim_{k \rightarrow \infty} \left(\sup_{j\in\bZ} h( \Phi( k, \sigma^{j}(\bu), B[j] ) ,  x[j + k ]  )\right) = 0
    \end{equation}
    then we say $\{ B[k] \}_{k \in \bZ}$ is \emph{uniformly attracted} to $ \{ x[k] \}_{k \in \bZ}$.
    \item[(ii)]
    We say $ \{ x[k] \}_{k \in \bZ}$ is a UAES if there is a neighbourhood $\{ B[k] \}_{k \in \bZ}$ of $ \{ x[k] \}_{k \in \bZ}$ that is uniformly attracted to $ \{ x[k] \}_{k \in \bZ}$.
\end{itemize}
\end{defn}
\begin{remark}
Note that if $ \{ x[k] \}_{k \in \bZ} $ is a UAES, then it is both a forward attractor and a pullback attractor \cite[Definition 3.11]{kloeden2011nonautonomous}.
Moreover, note that a neighbourhood $\{ B[k] \}_{k \in \bZ}$ of $ \{ x[k] \}_{k \in \bZ}$ has necessarily positive Lebesgue measure.
Therefore a UAES $ \{ x[k] \}_{k \in \bZ}$ cannot be an entire solution attracting a nonautonomous set which has some fibres with zero Lebesgue measure.
\end{remark}

UAESs allow us to rigorously define the number of stable RNN responses as a function of the specific input sequence driving the dynamics.
\begin{defn}
\label{def:decomp}
We say the system \eqref{eq:rnn_nds} with input $\bu$ admits a \emph{decomposition} into $n\geq 1$ UAESs if there are $n$ UAESs $\{ x_1[k] \}_{k \in \bZ}, \ldots, \{ x_n[k] \}_{k \in \bZ} $ such that, for all $\eta>0$ and $i=1,\ldots,n$, there are neighbourhoods $ \{ B_i^\eta[k] \}_{k \in \bZ}$ uniformly attracted by $ \{  x_i[k] \}_{k \in \bZ}$ and
\begin{equation}
\label{eq:cover_measure}
    \lambda(X \setminus \bigcup_{i=1}^n B_i^\eta[k] ) < \eta,   \qquad \forall k\in \bZ,
\end{equation}
where $\lambda$ denotes the Lebesgue measure on $X\subset \bR^{N_r}$. We say this is a \emph{proper decomposition} if in addition
\begin{equation}
    \label{eq:proper_decomp}
    \inf_{k \in \bZ} d_X( x_i[k] , x_j[k] ) > 0, \qquad \forall i,j \in \{ 1, \ldots, n \}, ~i\neq j.
\end{equation}
\end{defn}
Note that \eqref{eq:uni_pt_att} together with \eqref{eq:proper_decomp} imply that $ B_i^\eta[k] \cap  B_j^\eta[k] = \emptyset  $ at each time step $ k \in \bZ $. 
Condition \eqref{eq:cover_measure} means we can exclude a neighbourhood of the repelling dynamics that may sit on basin boundaries, while \eqref{eq:proper_decomp} implies there is a uniform threshold that can be used to separate all attractors.
\begin{defn}
	\label{def:uniform_echo_index}
	We say the system \eqref{eq:rnn_nds} driven by input $\bu \in \cU$ has \emph{echo index} $n\geq 1$, and write 
	$$
	\cI({\bf u})=n,
	$$
	if it admits a proper decomposition into $n$ UAESs. In this case, we say \eqref{eq:rnn_nds} has the $n$\emph{-Echo State Property} ($n$-ESP) for input $\bu$.
\end{defn}
Note that if $\bu$ has echo index $n$, then so does $\sigma^m(\bu)$ for any $m\in\bZ$.
\begin{figure}
    \centering
    \includegraphics[keepaspectratio=true,scale=0.65]{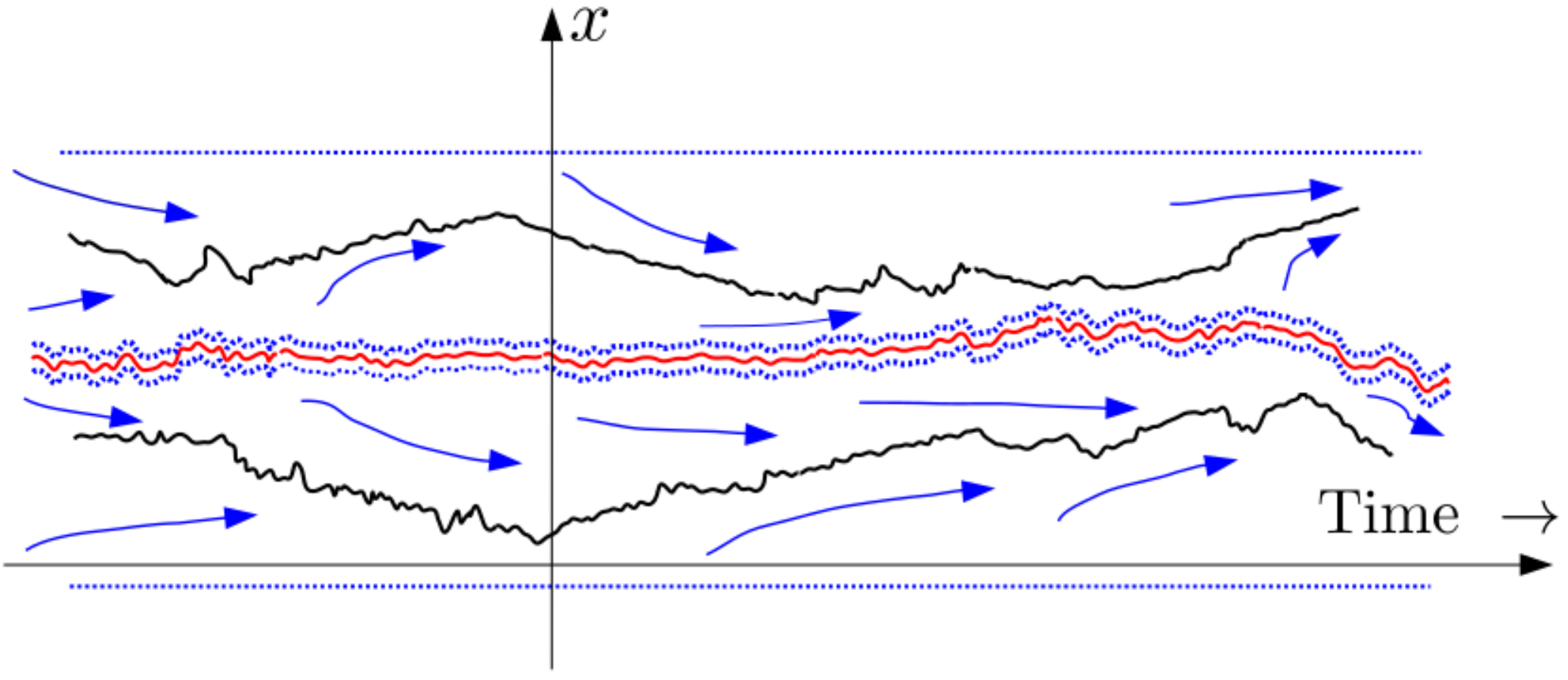}
    \caption{Blue dashed curves represent the boundaries of two nonautonomous compact sets which are uniformly attracted by two UAESs depicted as black curves. In red is depicted an entire solution that delimits the boundary between the basins of attraction of the two UAESs. Uniform convergence is guaranteed as long as, at each time step, we exclude a small neighbourhood of the separatrix solution shown in red;
    the smaller the excluded neighbourhood of the separatrix solution, the smaller the residual of the measure of \eqref{eq:cover_measure}.}
    \label{fig:echo_index}
\end{figure}

Definition \ref{def:uniform_echo_index} describes the case where a finite number $\cI(\bu)$ of ``attracting'' solutions (here modeled as UAESs) emerge as a consequence of driving the system with input sequence $\bu$.
Apart from a zero measure set of initial conditions, the RNN dynamics converge to one of these UAESs: which of these is actually followed will depend on the particular initial condition, see Figure \ref{fig:echo_index}.

Our definition of echo index assumes that all local nonautonomous attractors are UAESs. This is a restriction even when the input is null. In fact, there can be invariant curves or chaotic sets that attract some portion of the phase space. 
For all those cases where a proper decomposition in $n$ UAESs does not exist, we say the driven system \eqref{eq:rnn_nds} has \emph{indefinite} echo index.
We refer the reader to \cite{caraballo2013morse} for phase space decompositions in more general nonautonomous attractors, and to \cite{crauel2004towards} for decompositions framed within the theory of random dynamical systems.

Condition \eqref{eq:proper_decomp} guarantees that the UAESs do not merge into each other, implying that the number $n$ of stable responses to an input in the infinite past coincides with the number of stable responses in the infinite future. Hence, the echo index characterises those input-driven systems whose degree of multistability is invariant over time.
\begin{remark}
The set of input-driven systems possessing the ESP according to Definition \ref{def:uniqueentiresolution-ESP} and those having echo index 1 do not coincide. 
The example in Section \ref{sec:kloeden} represents an input-driven system possessing the ESP according to Definition \ref{def:uniqueentiresolution-ESP} whilst, however, a UAES does not exist. Therefore the echo index is not well-defined for such cases.
Vice-versa, there might be cases where input-driven systems have one UAES attracting a nonautonomous set with full Lebesgue measure (case with echo index 1), but such a UAES might coexist with other entire solutions attracting some zero measure set of phase space. In such a scenario, there will be more than one entire solution and thus Definition \ref{def:uniqueentiresolution-ESP} is not satisfied. Nevertheless, such an input-driven system will produce a unique stable response for almost all initial conditions.
\end{remark}

It is possible to extend the definition of echo index to sets of input sequences as follows.
\begin{defn}
\label{def:index-set}
Consider a system \eqref{eq:rnn_nds} with a set of possible input sequences $\cV\subset \cU$ and define
\begin{equation}
\label{eq:class_index}
    \cV_{[n]}=\{{\bf v}\in \cV~:~\cI({\bf v})=n\}.
\end{equation}
We say that such a system has echo index set $\cN(\cV)$ for a subset $\cV\subset \cU$ if
\begin{equation}
    \label{eq:index-set}
    \cN(\cV)=\{n\in\bN~:~ \cV_{[n]}\neq \emptyset \}.
\end{equation}
\end{defn}

We can split a set of inputs $\cV$ into a disjoint union according to the echo index:
\begin{equation}
    \label{eq:index-decomposition}
    \cV=  \left(\bigcup_{n\in\cN(\cV)} \cV_{[n]}\right)\cup \cV_{[\text{ind}]},
\end{equation}
where $\cV_{[\text{ind}]}$ are those input sequences that give an indefinite echo index.

\section{Theoretical results}
\label{sec:theoretical_results}

\subsection{A theorem of existence and uniqueness for uniformly attracting entire solutions}
\label{eq:unique_local_global}

First, let us delimit the region of uniform contraction of the map $G: U \times X \longrightarrow X $.
\begin{defn}
\label{def:contracting_set}
Given a positive real number $\mu$ such that $ 0< \mu < 1 $, we define the set of \emph{linear $\mu$-contraction uniform in $U$} of the map $G$ as 
\begin{equation}
    \label{eq:contracting_set}
    C(\mu,U):= \{ x \in X \,\, : \,\, \sup_{u \in U} \lVert D_x G(u,x)\rVert\leq \mu    \},
\end{equation}
where $\Vert\cdot\rVert$ denotes the matrix norm induced by the Euclidean norm on $\bR^{N_r}$.
\end{defn}
\begin{remark}
The set in \eqref{eq:contracting_set} is the phase space region where the nonautonomous dynamics contract at each time step with a rate of at most $\mu$, i.e. where each autonomous map $G(u,\cdot):X \rightarrow X $, with $u \in U$, contracts with a rate of at most $\mu$.
\end{remark}

It is possible to extend the definition of a positively invariant set $B\subset X $ to all input sequences in $\cV \subseteq \cU$.
\begin{defn}
\label{def:posit_invar}
A nonempty compact subset $B \subseteq X $ is called \emph{positively invariant} for $\cV \subseteq \cU$ (or $\cV$\emph{-positively invariant}) if
$$
\Phi(k, \bu, B) \subseteq B, \qquad \forall \bu \in \cV,~k\in\bZ_0^+
$$
\end{defn}
\begin{remark}
If the family of input sequences under consideration is $ \cV := \{ \sigma^{n}(\bv) \}_{n \in \bZ}  $, for some $\bv \in \cU$, then we simply say that $B$ is a $\bv$\emph{-positively invariant} set. 
This coincides with the case in Definition \ref{def:posit_invar_family} of a constant nonautonomous set $ B_n \equiv B $, which is positively invariant for the input $\bv$.
\end{remark}

The following theorem gives sufficient conditions to prove the existence and uniqueness of a UAES in a compact and convex phase space subset. The proof of the theorem exploits the contraction of the family of autonomous maps $ \{ G(u,\cdot) \}_{u \in U}$ defining the nonautonomous system. This result can be seen as a deterministic input-driven fixed point theorem, see \cite{smart1980fixed} for similar results in the context of autonomous dynamical systems or \cite{itoh1979random} in the context of random dynamical systems.
\begin{thm}
    \label{thm:idea_contraction}
    Let $\mu$ be a positive real number such that $ 0 < \mu < 1 $.
    Suppose $Q_\mu$ is a $\cU$-positively invariant nonempty compact set such that it is contained inside $ C(\mu,U) $ of Definition \ref{def:contracting_set}, and suppose further that $Q_\mu$ is convex.
    Then, for all $\bv \in \cU$ the system \eqref{eq:rnn_nds} driven by the input sequence $\bv$ admits a unique entire solution in $Q_\mu$. In particular, if such entire solution is contained inside the interior of $ Q_\mu $, then this entire solution is a UAES.
\end{thm}
\proof
The assumption of convexity of $ Q_\mu $ implies that for all $x_0, y_0 \in Q_\mu $ the line segment $ \ell_{[x_0,y_0]} := \{ z = (1-s)x_0+ s y_0 \,\, | \,\, s \in [0,1] \} $ lies inside $ Q_\mu $ and the following inequality holds
\begin{align}
    \label{eq:contractivity_argument}
    \begin{split}
    d_X( G(u,x_0 ) , G(u, y_0) ) =&  \lVert G(u,x_0 ) - G(u, y_0) \rVert   =  \left\lVert \int_{s=0}^{1} D_x G(u,(1-s)x_0+ s y_0) \cdot y \,\, ds \right\rVert \\
    \leq & \sup_{z \in \ell_{[x_0,y_0]} } \lVert D_x G(u, z ) \rVert \,\, \lVert x_0 - y_0 \rVert
    =  \sup_{z \in \ell_{[x_0,y_0]} }\lVert D_x G(u, z ) \rVert \,\, d_X( x_0, y_0 ) \leq \mu \,\, d_X( x_0, y_0 ) ,
    \end{split}
\end{align}
where the last inequality holds whenever $u \in U$, by the hypothesis that $ Q_\mu \subseteq C(\mu, U) $.

Now, let us take an input sequence $\bv \in U^\bZ$ and an arbitrarily chosen initial time step $k_0 \in \bZ $.
Denote $ x_k := \Phi (k, \sigma^{k_0}(\bv), x_0  ) $ and $ y_k := \Phi (k, \sigma^{k_0}(\bv), y_0  ) $, for some $k>0$.
Now, $\bv = \{ v[k] \}_{k \in \bZ}$ assumes values in $ U $, hence in particular $v[k_0] \in U $ and \eqref{eq:contractivity_argument} reads $  d_X( x_1 , y_1 ) \leq \mu \,\, d_X( x_0, y_0 )  $.
By hypothesis $Q_\mu$ is positively invariant for all input sequences $\bu \in U^\bZ$, thus $x_1, y_1 \in Q_\mu $ if $ x_0, y_0 \in Q_\mu $.
Therefore, we can repeat the same argument on the pair $(x_1,y_1)$ with $v[k_0+1]$ in place of $(x_0,y_0)$ with $v[k_0]$, obtaining that $  d_X( x_2 , y_2 ) \leq \mu \,\, d_X( x_1, y_1 ) \leq \mu^2 \,\, d_X( x_0, y_0 ) $, and so on for all pairs $( x_k, y_k) $ forming the forward trajectories.
In other words, by induction we proved that 
\begin{equation}
    \label{eq:cocycle_contraction}
  d_X( \Phi(k, \sigma^{k_0}(\bu),x_0), \Phi(k, \sigma^{k_0}(\bu), y_0) ) \leq \mu^k d_X( x_0, y_0 )\quad   \forall \bu \in \cU, \forall k_0 \in \bZ , \forall x_0,y_0 \in Q_\mu, \forall k >0.
\end{equation}

Moreover, for any given input sequence $\bv \in \cU $, thanks to the fact that $Q_\mu $ is $\bv$-positively invariant, the nonautonomous set defined as 
\begin{equation}
\label{eq:pullback_entire}
    A_n  := \bigcap_{m \leq n} \Phi( n-m , \sigma^{m}(\bv), Q_\mu ), \qquad  n \in \bZ 
\end{equation}
is an invariant nonautonomous set for input $\bv$ as in Definition \ref{def:invariant_nonaut_family} which is confined inside $Q_\mu$; see \cite[Lemma 2.20]{kloeden2011nonautonomous} for a proof of this fact.

We will prove that the invariant nonautonomous set of \eqref{eq:pullback_entire} is indeed a UAES.
Let us first prove that it is an entire solution.
By contradiction, let us assume there are $ x_s, y_s \in A_s $ distinct points, i.e. $ x_s \neq y_s $, for some $ s \in \bZ $.
By the invariance of $ \{ A_n \}_{n \in \bZ} $, it follows that $ \Phi( k, \sigma^{-k+s}(\bv), A_{-k+s}  ) = A_s $. Therefore, there exist $ x', y' \in A_{-k+s} $ such that $ x_s = \Phi( k, \sigma^{-k+s}(\bv), x'  ) $ and $ y_s = \Phi( k, \sigma^{-k+s}(\bv), y'  ) $.
Hence, the following inequalities hold:
$$
    0 < d_X ( x_s, y_s  ) =   d_X ( \Phi( k, \sigma^{-k+s}(\bv), x'  ), \Phi( k, \sigma^{-k+s}(\bv), y'  )  ) 
    \stackrel{\eqref{eq:cocycle_contraction}}{\leq} \mu^k d_X ( x', y' ) \leq \mu^k \diam(Q_\mu). 
$$
Therefore, by the compactness of $ Q_\mu $, there cannot exist two distinct points in $A_s$, i.e. the invariant nonautonomous set \eqref{eq:pullback_entire} is an entire solution $ \mathbf{x} = \{ x[n] \}_{n \in \bZ} $ for input $\bv$.
This fact implies that inside $Q_\mu$ there can be only one entire solution.
Indeed, if there exists another entire solution $\mathbf{y} = \{ y[n] \}_{n \in \bZ}  $ for input $\bv$, then it means that $ \Phi( k, \sigma^{-k+s}(\bv), y[-k+s]  ) = y[s] $ holds for all $s \in \bZ$ and for all $k > 0$.
Therefore, if such entire solution is contained in $ Q_\mu $, i.e. $y[n] \in Q_\mu$ for all $n \in \bZ$, then $ y[s] \in  \bigcap_{k \geq 0} \Phi( k , \sigma^{-k+s}(\bv), Q_\mu ) =   \bigcap_{m \leq s} \Phi( s-m , \sigma^{m}(\bv), Q_\mu ) = A_{s} $, which coincides with $ x[s]$. 

Now, we prove that such unique entire solution is uniformly attracting according to \eqref{eq:uni_pt_att} of Definition \ref{def:local_point_attractor}.
Let us take a sequence of integers $\{ j_k \}_{k \in \bZ} $.
Note that 
\begin{align}
      h( \Phi( k, \sigma^{j_k}(\bu), Q_\mu) ,  x[j_k + k ]  )  & =  \max_{x \in Q_\mu} d_X( \Phi( k, \sigma^{j_k}(\bu), x) ,  x[j_k + k ]  ) \\
                & =  \max_{x \in Q_\mu} d_X(  \Phi ( k, \sigma^{j_k}(\bv), x ) ,  \Phi ( k, \sigma^{j_k}(\bv), x[j_k] ) ) .
\end{align}
Now, relation \eqref{eq:cocycle_contraction} implies that $ \max_{x \in Q_\mu} d_X(  \Phi ( k, \sigma^{j_k}(\bv), x ) ,  \Phi ( k, \sigma^{j_k}(\bv), x[j_k] ) ) \leq \mu^k \diam( Q_\mu ) $ and hence that $\lim_{k \rightarrow \infty} h( \Phi( k, \sigma^{j_k}(\bu), Q_\mu) ,  x[j_k + k ]  ) = 0$, regardless of the sequence of integers $\{ j_k \}_{k \in \bZ} $. Therefore, $\lim_{k \rightarrow \infty} \sup_{n \in \bZ} h( \Phi( k, \sigma^{n}(\bu), Q_\mu) ,  x[n + k ]  ) = 0$.
In particular, all forward trajectories starting at the same time step from initial conditions inside $Q_\mu$ synchronise to a (unique) common solution which is indeed $\bx$.
Moreover, as long as $ x[n] \in  \interior( Q_\mu ) , \forall n \in \bZ $, then we can find a neighbourhood of $ \bx = \{ x[n] \}_{n \in \bZ}$ such that Definition \ref{def:local_point_attractor} is satisfied. 
\qed

\subsection{$n$-ESP for systems perturbed by low-amplitude inputs}
\label{sec:small_input}

In this section, we analyse the input-driven dynamics of an autonomous RNN perturbed by low-amplitude input sequences. A special case of \eqref{eq:rnn_nds} is the constant-input case:
\begin{equation}
x[k+1]=G( u_0, x[k]),
\label{eq:rnn_ndszero}
\end{equation}
i.e. with $ \bu \equiv u_0 $ constant, which gives an autonomous nonlinear system.
We will look at how the autonomous system \eqref{eq:rnn_ndszero}, possessing a uniformly stable point $x^*$, reacts under small external perturbations represented by a forcing input sequence.
Therefore, input values are taken in a neighbourhood of $ u_0 $, i.e. $\cU := B_r(u_0)^{\bZ}$ with $r>0$ but small.
The idea is that stable fixed points will become UAESs in the nonautonomous setting, giving rise to the $n$-ESP.

In the following, we will assume that the autonomous map \eqref{eq:rnn_ndszero} possesses a uniformly attracting stable fixed point according to the following definition.
\begin{defn}
\label{def:uasp}
Let us be given an autonomous map $F:X \longrightarrow X $ and a fixed point $x^* \in X$ for $F$, i.e. $F(x^*)=x^*$.
We call $x^*$ a Uniformly Attracting Stable Point (UASP) if there is an $0<M<1$ and $\delta >0 $, such that 
\begin{equation}
    \label{eq:uasp_F}
    d_{X}(F(z),x^*) < M d_{X}(z,x^*) \qquad \forall z \in B_\delta(x^*).
\end{equation}
\end{defn}

The property \eqref{eq:uasp_F} imposes some algebraic condition on the linearised map at the fixed point $x^*$.
\begin{lem}
\label{lem:max_sing_value}
Let us be given an autonomous map $F$ with one UASP $x^*$ characterised by a contraction rate $0<M<1$. 
Let $F$ be differentiable at $x^*$.
Denote by $ A :=  D_x F(x^*) $ the Jacobian matrix evaluated onto the fixed point $x^*$, and its maximum singular value as $\sigma(A)$.
Then $\sigma(A)\leq M $.
\end{lem}
A proof of Lemma \ref{lem:max_sing_value} can be found at the end of Appendix \ref{proof:lemma}.

\begin{thm}
\label{thm:small_input}
Suppose there is a UASP $x^*$ of \eqref{eq:rnn_ndszero}, characterised by a rate of contraction $M$ in a ball of radius $\delta$. Then, 
$ \forall \, 0<\varepsilon<1-M, \,\, \exists \,\, 0<\delta_\varepsilon < \delta, \, r_\varepsilon > 0 $ such that for all input sequences $\bu \in \overline{B_{r_\varepsilon }(u_0)}^\bZ$ there exists a unique UAES which uniformly attracts $B_{\delta_\varepsilon}(x^*)$ with a rate of $M+\varepsilon$.
\end{thm}

\proof
We claim that $ \forall \, 0<\varepsilon<1-M, \,\, \exists \, 0<\delta_\varepsilon<\delta \,\,:\,\, \forall \, 0< \gamma < \delta_\varepsilon(1-M), \,\, \exists \, \zeta_{\gamma, \varepsilon} > 0 $ such that for all input sequences $\bu \in \overline{B_{\zeta_{\gamma, \varepsilon}}(u_0)}^\bZ$ there exists a unique UAES confined in $\overline{B_{\frac{\gamma}{1-M}}(x^*)}$ which uniformly attracts the neighbourhood $B_{\delta_\varepsilon}(x^*)$ with a rate of $M+\varepsilon$. 
This will prove the thesis choosing $r_\varepsilon$ in the statement of Theorem \ref{thm:small_input} as $ \zeta_{\gamma, \varepsilon} $.\\
By definition, $\lVert A\rVert=\sigma(A)$, where $\lVert\cdot\rVert$ denotes the matrix norm induced by the Euclidean norm on $\bR^{N_r}$.
Lemma \ref{lem:max_sing_value} and the fact that $x^*$ is a UASP for $G(u_0,\cdot):X\rightarrow X$ imply that $ \|D_x G(u_0,x^*)\|\leq M$.
Assumption \ref{assum:fundamental}(i) implies continuity of $\|D_x G\| $ at $(u_0,x^*)$, hence for all $\varepsilon>0$ we can find a $ \rho_\varepsilon > 0  $ and a $\delta_\varepsilon>0$ such that $ \max_{u' \in \overline{ B_{\rho_\varepsilon}(u_0) } } \lVert D_x G(u', x ) \rVert \leq M + \varepsilon ,  \quad  \forall x \in \overline{B_{\delta_\varepsilon}(x^*) } $.
In terms of Definition \ref{def:contracting_set} we can state that for all $0<\varepsilon<1-M$ and denoting $\mu=M+\varepsilon <1$, there exist $\rho_\varepsilon, \delta_\varepsilon >0$ such that $ \overline{B_{\delta_\varepsilon}(x^*) } \subseteq C(\mu, \overline{ B_{\rho_\varepsilon}(u_0) } ) $.
Now, continuity of $G(\cdot, x) : U \rightarrow X $ for each $x \in \overline{B_{\delta_\varepsilon}(x^*) }  $ implies that for all $\gamma>0 , \exists \zeta_\gamma>0$ such that 
\begin{equation}
    \label{eq:G_continuous}
    G( \overline{B_{\zeta_\gamma}(u_0) } , \overline{B_{\delta_\varepsilon}(x^*) }  ) \subseteq B_\gamma ( G(u_0, \overline{B_{\delta_\varepsilon}(x^*) ) } .
\end{equation}
The hypothesis that $x^*$ is a UASP for the autonomous map $F(\cdot)=G(u_0, \cdot)$ reads $ G(u_0 , \overline{B_{\delta_\varepsilon}(x^*) } )  \subseteq  B_{M\delta_\varepsilon}(x^*) $, as long as $\delta_\varepsilon < \delta$. Therefore, we obtain from \eqref{eq:G_continuous} that 
\begin{equation}
    \label{eq:above_argument}
    G( \overline{B_{\zeta_\gamma}(u_0) } , \overline{B_{\delta_\varepsilon}(x^*) }  ) \subseteq   \overline{B_{M \delta_\varepsilon + \gamma }(x^*)  } .
\end{equation}
Therefore, as long as we choose $\gamma $ such that $ M\delta_\varepsilon + \gamma < \delta_\varepsilon $, i.e. $ \gamma < \delta_\varepsilon(1-M) $, we get from \eqref{eq:above_argument} that $ G( \overline{B_{\zeta_\gamma}(u_0) } , \overline{B_{\delta_\varepsilon}(x^*) }  ) \subset  \overline{B_{\delta_\varepsilon}(x^*) }  $, which proves that $ \overline{B_{\delta_\varepsilon}(x^*) } $ is a $\cU$-positively invariant set, considering $\cU=\overline{B_{\zeta_\gamma}(u_0) }^\bZ$.
This allows us to make use of Theorem \ref{thm:idea_contraction} with $ Q_\mu= \overline{B_{\delta_\varepsilon}(x^*)}$ and $\cU=\overline{B_{\zeta_{\gamma, \varepsilon}}(u_0)}^\bZ$, with $\zeta_{\gamma, \varepsilon}=\min\{ \rho_\varepsilon, \zeta_\gamma \}$.
Hence, we proved that $ \forall \, 0<\varepsilon<1-M, \,\, \exists \, 0<\delta_\varepsilon<\delta \,\,:\,\, \forall \, 0< \gamma < \delta_\varepsilon(1-M), \,\, \exists \, \zeta_{\gamma, \varepsilon} > 0 $ such that for all input sequences $\bu \in \overline{B_{\zeta_{\gamma, \varepsilon}}(u_0)}^\bZ$ there exists a unique entire solution $\bx$ confined in $ \overline{B_{\delta_\varepsilon}(x^*) } $.
Moreover, note that once chosen a $\gamma $ such that $ M\delta_\varepsilon + \gamma < \delta_\varepsilon $ we can iterate the above argument in \eqref{eq:above_argument} leading to $ G( \overline{B_{\zeta_\gamma}(u_0) } , \overline{B_{M \delta_\varepsilon + \gamma }(x^*) }  ) \subseteq   \overline{B_{M^2 \delta_\varepsilon + \gamma(1+M) }(x^*)  }  $ and so on.
This leads to a sequence $ \delta_{n} = M \delta_{n-1} + \gamma $, for all $ n \geq 1$, starting with $ \delta_0 = \delta_\varepsilon $, which is strictly decreasing. Unfolding such a recursive relation leads to $\delta_n = M^n\delta_\varepsilon + \gamma\sum_{j=0}^{n-1}M^j$, whose limit is $\dfrac{\gamma}{1-M}$.
As a consequence $\bx$ is confined in $\overline{ B_{\frac{\gamma}{1-M}}(x^*) } $ which is strictly contained in the interior of $ \overline{B_{\delta_\varepsilon}(x^*) } $; thus $\bx $ is an UAES uniformly attracting the neighbourhood $ B_{\delta_\varepsilon}(x^*) $ with a rate of $M+\varepsilon$.
\qed

The assumption on $x^*$ to be a UASP makes proof of existence of a UAES easier.
Nevertheless, the same results can presumably be proved just assuming $x^*$ to be a stable hyperbolic fixed point, using exponential dichotomies; see for example \cite{potzsche2011nonautonomous}.

As a corollary of Theorem \ref{thm:small_input} we get the following result.
\begin{corol}
 If for a given constant input value $u_0$ the autonomous map \eqref{eq:rnn_ndszero} presents a number $n\geq1$ of UASPs, then the nonautonomous system $ x[k] = G(u[k], x[k-1]) $ driven by any deterministic input sequence $u[k]$ assuming values in a neighbourhood of $u_0$ will present (at least) $n$ UAESs, as in Definition \ref{def:local_point_attractor}.
 Therefore, input sequences $\bu \in B_r(u_0)^\bZ$ with small enough $r>0$ will typically give rise to an input-driven system with the $n$-ESP.
\end{corol}

\subsection{ESP for RNNs driven by large-amplitude inputs}
\label{sec:esp-large-inputs}

We rigorously prove, in Proposition~\ref{prop:large_inputs}, that driving RNNs with large-amplitude inputs generally leads to echo index 1. Our result is similar to the one obtained in \cite[Theorem 2]{manjunath2013echo}, with the important differences that we (i) account for the uniform (i.e. also forward in time) attractiveness of solutions, (ii) consider also the feedback of the output, and (iii) highlight how $W_{in}$ induces a geometric structure in the space of input values relevant for the analysis of reliable responses of an input-driven RNN.

First, we prove a technical Lemma that requires some additional definitions.
We will denote with $ (W)_{(j)} $ the $j$th row of a matrix $W$.
\begin{defn}
We define
\begin{equation}
    \label{eq:hyperplane}
    H_j:= \{ u \in \mathbb{R}^{N_i} \,\, : \,\, (W_{in})_{(j)} \cdot u = 0 \}
\end{equation}
the hyperplane of input values which vanish the $j$th row of the matrix $W_{in}$.
Moreover, given real values $\epsilon>0$ (meant to be small) and $R>0$ (meant to be large), we define
\begin{equation}
    \label{eq:ESP_region}
    P_j(\epsilon, R ) := \left\{ u \in \bR^{N_i}\setminus B_R(0) \,\, : \,\,  \dfrac{ \left| (W_{in})_{(j)} \cdot u \right| }{ \lVert (W_{in})_{(j)} \rVert \lVert u \rVert }   \geq \epsilon   \right\} .
 \end{equation}
\end{defn}
\begin{remark}
\label{rem:interpretation}
Note that, denoted with $\theta$ the angle between the vectors pointing to $(W_{in})_{(j)}$ and $ u $, then $ \cos(\theta) = \dfrac{(W_{in})_{(j)} \cdot u}{\lVert (W_{in})_{(j)} \rVert \lVert u \rVert}   $.
Therefore, the set $P_j(\epsilon, R )$ is basically the set $ \bR^{N_i}\setminus B_R(0)  $ where we cut out all the lines which form an angle $ \theta \in \left(\arccos(\epsilon), \arccos(-\epsilon) \right)  $, i.e. those close to the hyperplane $H_j$.
\end{remark}
\begin{lem}
\label{lem:xi_functions}
Let $W_r, W_{fb}, W_{in} $ be real matrices of dimensions, respectively, $ N_r\times N_r, \,  N_r\times N_o, \, N_r \times N_i $, and $\psi \in C^0(\mathbb{R}^{N_r},\mathbb{R}^{N_o})$.
Consider the functions $ \xi_j: \mathbb{R}^{N_i} \times [-L,L]^{N_r} \longrightarrow \mathbb{R} $, for $j=1,\ldots, N_r$, defined as $ \xi_j (u,x) := (W_{in})_{(j)} \cdot u +  f_j(x) $, where $ f_j(x) := (W_r)_{(j)}\cdot x + (W_{fb})_{(j)} \cdot \psi(x) $.
If $(W_{in})_{(j)}$ is not the null vector then in the subset of input values $ P_j(\epsilon, R ) $ we can make the function $| \xi_{j}(u,x) |$ large as much as we want. Precisely, for all $\epsilon>0$ and for all $\overline{\xi}>0$ there exists an $ R_{\overline{\xi},\epsilon}>0 $ (which depends on $ \| (W_{in})_{(j)} \|^{-1}$ and $\max_{x \in [-L,L]^{N_r}} |f_j(x)|   $) such that for all $R\geq R_{\overline{\xi},\epsilon}$ we have that $ \inf \left\{ \,\,| \xi_{j}(u,x) | \,\, : \,\, u \in P_j(\epsilon, R ) \right\} \geq \overline{\xi} $ holds for all $ x \in [-L,L]^{N_r} $.
\end{lem}

\proof
If $(W_{in})_{(j)}$ is the null vector then $\xi_{j}(u,x)=f_j(x)$, i.e. it does not depend on the input, thus let us assume $(W_{in})_{(j)}$ is not the null vector.
Then, a hyperplane $H_j= \{ u \in \mathbb{R}^{N_i} \,\, : \,\, (W_{in})_{(j)} \cdot u = 0 \}$ is defined in the space of input values $\mathbb{R}^{N_i}$ such that $H_j$ is the orthogonal space of the vector pointing to $ (W_{in})_{(j)} $.

In general, 
$$ 
    (W_{in})_{(j)} \cdot u = \lVert (W_{in})_{(j)} \rVert  R \cos(\theta) 
$$
where $R$ is the norm of $u$ and $\theta$ is the angle between the vectors pointing to $(W_{in})_{(j)}$ and $ u $. Let us fix arbitrarily a $R>0$ and consider $u$ to vary on the surface of the ball of radius $R$ centred on the origin of $\bR^{N_i}$, i.e. $ u \in \partial B_R(0)$, thus we have that $ \| (W_{in})_{(j)} \| R $ is constant and $\cos(\theta) $ can assume any value in $[-1,1]$. 
In particular, $ (W_{in})_{(j)} \cdot u $ can assume any value in the interval $ \left[-\| (W_{in})_{(j)} \| R, \,\, \| (W_{in})_{(j)} \| R \right] $.

Now, let us take an arbitrary $\epsilon>0$ and consider the subset $ P_j(\epsilon, R ) $ for some $R>0$.
Note that for any given $ u \in P_j(\epsilon, R )$, denoted with $\varepsilon := \left| \cos{ \left( \theta \right)} \right| $, where $\theta $ is the angle between $  (W_{in})_{(j)} $ and $u$, and with $\rho$ the norm of $\|u\|$, then $\varepsilon\geq \epsilon$ and $\rho \geq R$.
Therefore, we have that 
\begin{equation}
    \label{eq:inf_bound}
    | (W_{in})_{(j)} \cdot u |  \geq \| (W_{in})_{(j)}\| R \epsilon 
\end{equation}
holds for all $u \in P_j(\epsilon, R )$.
Therefore, the reverse triangle inequality leads to
\begin{equation}
\label{eq:xi_reverse}
    | \xi_{j} (u,x)| = | f_j(x)  + (W_{in})_{(j)} \cdot u   |   \geq \Bigl|  \bigl| f_j(x) \bigl| - \bigl| (W_{in})_{(j)} \cdot u \bigl|  \Bigl|  
                        \geq \bigl| (W_{in})_{(j)} \cdot u \bigl| \,\, -  \bigl| f_j(x) \bigl|,
\end{equation}
where the last inequality holds for all $x \in [-L,L]^{N_r}$ and for all $u \in P_j(\epsilon, R )$ as long as $R \geq \dfrac{\sigma}{\epsilon \| (W_{in})_{(j)} \|}$, where we denoted $\sigma=\max_{x \in [-L,L]^{N_r}} |f_j(x)|$. \\
Furthermore, for all $\overline{\xi}>0$ if $R \geq R_{\overline{\xi}, \epsilon}:=  \dfrac{ \overline{\xi} + \sigma }{\epsilon \| (W_{in})_{(j)} \|} $ then thanks to \eqref{eq:inf_bound} and \eqref{eq:xi_reverse} we get that the following
$$
    | \xi_{j} (u,x)| \geq \bigl| (W_{in})_{(j)} \cdot u \bigl| \,\, -  \bigl| f_j(x) \bigl| \geq \| (W_{in})_{(j)} \| R_{\overline{\xi}, \epsilon} \epsilon  -\sigma \geq \overline{\xi}
$$
holds for all $x \in [-L,L]^{N_r}$ and for all $ u \in P_j(\epsilon,R)$.
In conclusion, we have that for all $\epsilon>0$ and for all $\overline{\xi}>0$ there exists an $ R_{\overline{\xi},\epsilon}>0 $ such that for all $R \geq R_{\overline{\xi},\epsilon}$ we have that
$$
    \inf_{u \in P_j(\epsilon, R )  } | \xi_{j}(u,x) |  \geq \overline{\xi}
$$
holds for all $ x \in [-L,L]^{N_r} $.
\qed

We now state the main result of this section.
\begin{prop}
\label{prop:large_inputs}
Consider the RNN \eqref{eq:leaky_rnn}-\eqref{eq:output} with a bounded $\phi\in C^1(\bR,(-L,L))$ that is monotonically increasing and $\phi'$ has a unique maximum point at $\xi=0$ and $\psi \in C^1(\mathbb{R}^{N_r}, \bR^{N_o})$. 
Then, the following two statements are true:
\begin{itemize}
    \item[(i)] If the following condition holds, for some $0<\mu<1$,
        \begin{equation}
        \label{eq:ESP_max_sing_value}
            \phi' ( 0 )  \lVert W_r + W_{fb}  D_x \psi(x) \rVert \leq \mu, \qquad \forall x \in [-L,L]^{N_r},
        \end{equation}
        then for all sets of input values $ V \subseteq \mathbb{R}^{N_i}$ it holds that $ X \subseteq C(\mu,V)$ of \eqref{eq:contracting_set}. 
        In particular, Theorem \ref{thm:idea_contraction} implies that, for all compact sets $U \subset \bR^{N_i}$ and for all input sequences $\bu \in U^\bZ$, the RNN \eqref{eq:leaky_rnn}-\eqref{eq:output} driven by $\bu $ admits exactly one UAES which attracts the whole phase space $X=[-L,L]^{N_r}$, i.e. it has echo index 1.
        
    \item[(ii)] If each row $ (W_{in})_{(j)} $ is non zero, 
    then $\forall \varepsilon >0 $ and $0<\mu<1$ there exists a radius $R_{\varepsilon, \mu} > 0 $ such that $X\subseteq C(\mu,V) $ as long as the dynamics are driven by inputs assuming values in a set $V $ such that $V \subseteq  \bigcap_{j=1}^{N_r} P_j(\epsilon, R_{\epsilon,\mu}) $ of \eqref{eq:ESP_region}.
    In particular, for any compact $U \subset \bigcap_{j=1}^{N_r} P_j(\epsilon, R_{\epsilon,\mu})$, Theorem \ref{thm:idea_contraction} implies that for any $\bu \in U^\bZ$ there exists a unique UAES which attracts the whole phase space $X$, i.e. it has echo index 1. 
\end{itemize}
\end{prop}
\proof
First, let us consider the case of $\alpha=1$.
Given a subset of input values $ V\subseteq \bR^{N_i}$, we define
$$
    l_V(x) : = \sup_{u \in V} \lVert S(u,x) \rVert,
$$
where $S(u,x)$ is given in \eqref{eq:diag_deriv}.
Then, for all $x\in [-L,L]^{N_r}$, the norm of $D_x G$ \eqref{eq:jacobian_G}, can be upper-bounded with:
$$
   \sup_{u \in V} \lVert D_x G(u,x) \rVert  \leq l_V(x) \, \lVert M(x) \rVert .
$$
Now, since $ S(u,x) $ is a diagonal matrix, we have
\begin{align}
    l_V(x) = \sup_{u \in V} \| S(u,x)\|
        & = \sup_{u \in V} \max_{j=1,\ldots,N_r} | \phi'( \xi_j(u,x) ) | \\
        & = \sup_{u \in V} \phi' \left( \min_{j=1,\ldots,N_r} |\xi_j(u,x)| \right)\\
        & = \phi'\left(  \inf_{u \in V} \left\{ \min_{j=1,\ldots,N_r} \bigl| \xi_j(u,x)\bigl| \right\} \right) ,
\end{align}
where the last two equalities hold because $\phi'$ is a continuous function with a unique maximum point in $\xi=0$.
\begin{itemize}
    \item[(i)] Now, since $\phi'(z)$ assumes its maximum value at $z=0$ then $l_V(x)\leq \phi'(0)$. Therefore, if there exists a $0<\mu<1$ such that \eqref{eq:ESP_max_sing_value} holds then 
    $$
       \sup_{u \in V} \lVert D_x G(u,x) \rVert  \leq l_V(x) \, \lVert M(x) \rVert  \leq \phi' ( 0 )  \lVert M(x) \rVert \stackrel{\eqref{eq:ESP_max_sing_value}}{\leq} \mu .
    $$
    In other words, if \eqref{eq:ESP_max_sing_value} holds then $X \subseteq C(\mu, V) $ for any $V \subseteq \bR^{N_i}$.
    Now, since $X=[-L,L]^{N_r}$ is a convex and compact set which is $ U^\bZ$-positively invariant for any compact $U \subset \bR^{N_i}$ (see Proposition \ref{prop:absorbing-set-pullback}), we conclude the claim, applying Theorem \ref{thm:idea_contraction} with $\cU:= U^\bZ $ on the whole space $X=[-L,L]^{N_r}$.
    
    \item[(ii)] Let us fix an arbitrary $\epsilon>0$ and a $\mu \in (0,1)$. We will prove that there exists a $ R_{\epsilon, \mu} >0 $ such that $\sup_{u \in V} \lVert D_x G(u,x) \rVert \leq \mu \text{ w.r.t. input values in a subset } V \subseteq \bigcap_{j=1}^{N_r} P_j(\epsilon, R_{\epsilon,\mu})  $.
    First of all, thanks to Lemma \ref{lem:xi_functions} we know that for all $\epsilon>0$ and $\overline{\xi}>0$ there exists an $ R_{\overline{\xi}, \epsilon} > 0 $ such that 
    $$
        \inf_{ u \in P_j(\epsilon, R_{\overline{\xi}, \epsilon}) } \Bigl\{ \min_{j=1,\ldots,N_r} \bigl| \xi_j(u,x)\bigl| \Bigl\} \geq \overline{\xi}, \qquad \forall x \in [-L,L]^{N_r},
    $$
    as long as each row $ (W_{in})_{(j)} $ is not the zero vector.
    Moreover, $ \phi \in C^1(\mathbb{R}, (-L,L) ) $ is monotonic and its image is $(-L,L)$, hence $ \lim_{|\xi|\rightarrow{\infty}} \phi'(\xi)=0$. 
    Hence, $ l_{V}(x)$ can be made arbitrarily low, regardless of $x$, as long as $ V \subseteq P_j(\epsilon, R)  $ with $R$ large enough. 
    More precisely, denoting $ \Tilde{\sigma}:=  \max_{x \in [-L,L]^{N_r}}\|M(x)\| $, for all $\epsilon>0$ and $0<\mu<1$ there exists an $R_{\epsilon, \mu}>0$, such that $ l_{V}(x) \leq \dfrac{\mu}{\Tilde{\sigma} } $ (an consequently such that $\sup_{u \in V} \|D_x G(u,x)\| \leq \mu$) holds for all $x \in  [-L,L]^{N_r} $,  w.r.t. input values in a subset $V \subseteq \bigcap_{j=1}^{N_r} P_j(\epsilon, R_{\epsilon,\mu})  $.
    Finally, Theorem \ref{thm:idea_contraction} implies that for any compact $ U\subset \bigcap_{j=1}^{N_r} P_j(\epsilon, R_{\epsilon,\mu})   $ the RNN driven by any $\bu \in U^\bZ$ admits a UAES which attracts the whole phase space, i.e. it has echo index 1.
\end{itemize}

On the other hand, when $\alpha\in(0, 1)$, we have
$$
    \|D_x G(u,x)\| \leq  (1-\alpha)^{N_r} +  \alpha^{N_r} \|S(u,x) M(x)\| \leq 1-\alpha  \bigl( 1 - \|S(u,x) M(x)\| \bigl),
$$
where the last inequality holds since the function $ f(x) = (1-\alpha)^x +  \alpha^x s $ is strictly monotonically decreasing for all $\alpha \in (0,1)$ and for all $ s \geq 0$.
Now, from these last inequalities we note that if $x \in [-L,L]^{N_r}$ is such that $  \sup_{u \in V} \lVert S(u,x) M(x) \rVert \leq \mu \in (0,1)$ then it holds that $ \sup_{u \in V} \lVert D_x G(u,x)\rVert\leq 1-\alpha(1-\mu) $, which leads back to the already proved case with $\alpha=1$.
Therefore, we can apply Theorem \ref{thm:idea_contraction} on the whole space $[-L,L]^{N_r}$ which is contained inside $C\bigl(1-\alpha(1-\mu), V \bigl)$.
\qed
\begin{remark}
Remarkably, the result of Proposition~\ref{prop:large_inputs}(ii) does not rely on any particular assumption on the matrices $W_r , W_{fb}$, and the readout $ \psi $, as long as the input values are large enough in amplitude and are sufficiently far from $ H=\bigcup_{j=1}^{N_r}  H_j$, in a sense made precise in Remark~\ref{rem:interpretation}.
The need to exclude the sets $ H_j $ is due to the fact that, for all input sequences $ \bu \in H_j^\bZ $ in some phase space directions, the RNN is basically a mere function of $x$.
Hence, in general, in these directions the dynamics might be repulsive unless we impose some conditions on the matrices $W_r, W_{fb}$ and the function $\psi$, as the condition in \eqref{eq:ESP_max_sing_value}.
Nevertheless, note that the set $ H = \bigcup_{j=1}^{N_r}H_j $ is the union of hyperplanes in $\mathbb{R}^{N_i}$, which has a zero Lebesgue measure in the space of input values.
Roughly speaking, this means that if the input sequence $\bu$ is a realisation of a random process which generates values inside $U$ according to a uniform distribution, then for all $R>0$ such that $ U \subset \bR^{N_i} \setminus B_R(0)  $ the probability to observe values of the input outside of the compact space $ U \cap \bigcap_{j=1}^{N_r}P_j(\epsilon,R)$ will be proportional to $\epsilon$.

Note that, in the presence of feedback of the output with linear readout, i.e. $ \psi(x) = W_o x $, whenever the activation function is such that $\phi'(0)=1$, as for example $\phi=\tanh$, then from \eqref{eq:ESP_max_sing_value} we obtain the condition on the maximum singular value $ \lVert M \rVert < 1 $, where $M$ is the effective recurrent matrix \eqref{eq:effective_recurrent_matrix}. In particular, when there is no feedback, i.e. $W_{fb}=0$, we recover the well-known sufficient condition $ \lVert W_r \rVert< 1 $.
\end{remark}

\subsection{Stability of echo index to inputs}
\label{sec:input-stable}

We now show how the metric imposed on the space of input sequences affects the stability of the echo index to input perturbations.
Notably, we show that the widely-used (e.g. \cite{grigoryeva2019differentiable}) product topology \eqref{eq:metric-inputsequences} implies that echo indices greater than 1 cannot remain the same (i.e. are not stable) even when considering arbitrarily small perturbations of input sequences.
\begin{thm}
\label{thm:dense} 
Let us consider an RNN of the form \eqref{eq:leaky_rnn}-\eqref{eq:output} driven by inputs in the space $\cU$ of bounded input sequences $\{ u[k] \}_{k \in \bZ}$, i.e. $\sup_k \{ d_{\bR^{N_i}}(0,u[k]) \} < +\infty$.
If $ (W_{in})_{(j)} $ is non zero for all $j = 1, \ldots, N_r$, then the subset $\cU_{[1]} \subset \cU$ \eqref{eq:class_index} of bounded input sequences for which an RNN driven by $\{ v[k] \}_{k \in \bZ}\in\cU_{[1]}$ has echo index 1 is dense in $\cU$ with the product topology \eqref{eq:metric-inputsequences}. 
\end{thm}

\proof
Assume an RNN of the form \eqref{eq:leaky_rnn}-\eqref{eq:output} satisfying hypothesis of Proposition \ref{prop:large_inputs}(ii) is given, i.e. such that  $ (W_{in})_{(j)} $ is non zero for all $j = 1, \ldots, N_r$.
Let's define $\cU:= \bigl\{ \{ u[k] \}_{k \in \bZ} \,\, : \,\, \sup_k \{ d_{\bR^{N_i}}(0,u[k]) \} < +\infty  \bigl\}$ as the universe of possible input sequences for driving the RNN.
Let us equip the space $\cU$ with the metric given by \eqref{eq:metric-inputsequences} and let's denote with $\cU_{[1]}$, as in \eqref{eq:class_index}, the subset of input sequences that gives echo index 1.

We claim that for all $\bu \in \cU$  and $\varepsilon>0$ there exists a $\bv \in \cU_{[1]}$ such that $d_{\text{prod}} ( \bu, \bv ) < \varepsilon $.

Let us take an arbitrary input sequence $\bu = \{ u[k] \}_{k \in \bZ}\in \cU$ and assume $\bu$ has not echo index 1. 
Proposition \ref{prop:large_inputs}(ii) says that $\forall \varepsilon >0 $ and $0<\mu<1$ there exists a radius $R_{\varepsilon, \mu} > 0 $ such that $X\subseteq C(\mu,U) $ as long as the dynamics are driven by inputs assuming values in a compact $U$ such that $U \subset  \bigcap_{j=1}^{N_r} P_j(\epsilon, R_{\epsilon,\mu}) $ of \eqref{eq:ESP_region}.
Now, let's define an input sequence $\bv = \{ v[k] \}_{k \in \bZ}$ such that $v[k]=u[k]$ for all $|k| \leq M$, for some $M>0$, and $v[k] \in U $ for all $ |k| > M $, for some compact $U \subset  \bigcap_{j=1}^{N_r} P_j(\epsilon, R_{\epsilon,\mu})  $ such that $U \subset B_{R_{\epsilon,\mu} + 1 }(0)$.
Now, driving the RNN with $\bv$ Proposition \ref{prop:large_inputs}(ii) implies that it is well defined a left-infinite internal state sequence $\mathbf{p}^{-M}:= \{ \ldots,p[-M-2], p[-M-1]\}$ which consists of the left branch (i.e. relatively to the infinite past) of the global pullback attractor for the input sequence $\bv$.
Therefore, at time-step $k= -M -1 $ we end up sitting exactly at the location $ p[-M-1] \in X $. 
Now, from the initial condition $x_0 = p[-M-1]$ we move driven by the sequence $\bv$ in the window of time $ [-M,M]$ tracing an orbit $ x[k] = \Phi(M+1+k, \sigma^{-M}(\bv),  x_0)$, for $k=-M, -M+1,\ldots, M$, in phase space and ending somewhere at position $ x[M]  \in X$.
Now, despite the potentially expanding dynamics in the window of time $[-M,M]$, Proposition \ref{prop:large_inputs}(ii) ensures that, when driven by $\bv$ for $k>M$ the entire phase space contracts forward in time with a rate $\mu$.
Hence, all forward trajectories starting from any initial condition at time step $k=M$ synchronise with the one starting from $x[M]$ in the future.
In particular, there exists a unique entire solution $\mathbf{z}$ for the system driven by $\bv$ and it is defined as
\begin{equation}
z[k] = 
 \begin{cases} 
      p[k] & k < -M \\
       \Phi(M+1+k, \sigma^{-M}(\bv),  p[-M-1]) & k \geq -M  .
   \end{cases}
\end{equation}

Moreover, $\mathbf{z}$ is a UAES which uniformly attracts the whole phase space.
Therefore, the RNN driven by $\bv$ has echo index 1. 
Finally, note that 
\begin{equation}
    \label{eq:}
        d_{\text{prod}}(\bu, \bv) =    \sum_{ |k|>M }\dfrac{d_U(u[k],v[k])}{2^{|k|}} \leq \diam(U)  \sum_{ |k|>M }\dfrac{1}{2^{|k|}}  \leq \diam(U) \dfrac{1}{2^{M+1}}   \leq    ( R_{\epsilon,\mu} + 1 )  \dfrac{1}{2^{M}}  
\end{equation}
can be made arbitrarily small by increasing $M$.
This proves that, for the topology induced by the metric $d_{\text{prod}}(\cdot, \cdot) $, the set $\cU_{[1]}$ is dense in $\cU$.
\qed  
\begin{remark}
As a consequence of Theorem \ref{thm:dense}, if an RNN has echo index $\mathcal{I}(\bu)\geq 2$ for input sequence $\bu$, then, considering the product topology induced by the metric in \eqref{eq:metric-inputsequences}, the echo index is not stable. 
In fact, the product topology allows arbitrarily large variations of input values in the far past and future of an input sequence $\bu$.
This means that there exists an arbitrarily small perturbation of $\bu$ leading to an input sequence $\Tilde{\bu}$ such that $\mathcal{I}(\Tilde{\bu})=1$. Interestingly, the proof of Theorem \ref{thm:dense} can be adapted to prove that the set $ \cU_{[\text{ind}]} $ is dense in $\cU$.
This observation suggests that the product topology is not suitable for the analysis of the stability of reliable responses of an RNN to variations of input sequences, e.g. due to noise, while the uniform metric \eqref{eq:uniform-metric-inputsequences} seems to be a more appropriate choice.
In fact, if a system \eqref{eq:rnn_nds} driven by an input sequence $\bu$ has a well-defined echo index, then its echo index might be stable when considering perturbations of $\bu$ which are uniformly bounded in time, as suggested by Theorem \ref{thm:small_input}.
\end{remark}

\section{Examples of input-driven RNNs with multistable dynamics}
\label{sec:experiments}

In this section, we consider the RNN in \eqref{eq:leaky_rnn} and a linear model for the RNN output \eqref{eq:output}.
The results here are illustrative of some of the behaviours we expect in more general cases.
In Section \ref{sec:2esp}, we show a two-dimensional example of RNN with echo index 2, highlighting the need for the concepts introduced in Definition~\ref{def:decomp}.
Then, in Section \ref{sec:Change-echo-index}, we show how the echo index may change also depending on the characteristics of the input sequence driving the dynamics. Section~\ref{sec:Maass} shows the application of our modeling framework to the high-dimensional task in \cite[Figure 5]{hoerzer2012emergence}.

\subsection{An example of switching system with echo index 2}
\label{sec:2esp}

We now report some numerical experiments on a simple RNN that can be thought of as switching dynamics between two autonomous maps.
We will consider the space of inputs defined as $\cU:=\{ u_1,u_2 \}^\bZ$, i.e. as made by sequences assuming just two possible values: $u_1$ and $u_2$.
Let us consider the map $G_{\alpha}(u,x):= (1-\alpha) x + \alpha \tanh(W_r x + W_{in} u) $ where 
$W_r = \begin{bmatrix}
    \frac{1}{2}  &  0 \\
    0     & \frac{3}{2}
\end{bmatrix}$,
$W_{in} = I_2$ is the identity matrix and $x$ a real vector of dimension two. What follows can be observed for any value of $\alpha \in (0,1]$. We select $\alpha=\frac{1}{4}$ because a small value slows down the state-update and hence highlights transient dynamics\footnote{As explained in equations (21)-(23) of \cite{ceni-ff-18}, $ \alpha $ scales the velocity field. In our example, the equation ruling the dynamics is $ x[k+1] = ( 1 -  \alpha  ) x[k] + \alpha \tanh( W_r x[k] + W_{in} u[k+1]  ) $ which can be equivalently written as 
$$
x[k+1] =  x[k] + \alpha \Bigl( \tanh( W_r x[k] + W_{in} u[k+1]  )  - x[k]  \Bigl) .
$$ 
In the latter expression is evident that the vector to be added to the current state $x[k]$ in order to get the next state $x[k+1]$ is scaled by the parameter $\alpha$. In this sense, a smaller $\alpha$ will slow down the dynamics highlighting the transients.}.

The space of possible input sequences is $\cU = \{ u_1, u_2 \}^{\bZ}$, where 
$ u_1:=  \begin{pmatrix}     
    \frac{1}{4}  \\
    \frac{3}{20}
\end{pmatrix} $ 
and $ u_2:=  - u_1 $. Therefore, the nonautonomous dynamics driven by input $\bu \in \cU$ consists of a switching pattern between the two component maps defined as $f_1(x):= G(u_1, x) $  and $f_2(x):= G(u_2, x) $.

\begin{figure}[ht!]
\begin{center}
\includegraphics[keepaspectratio=true,scale=0.45]{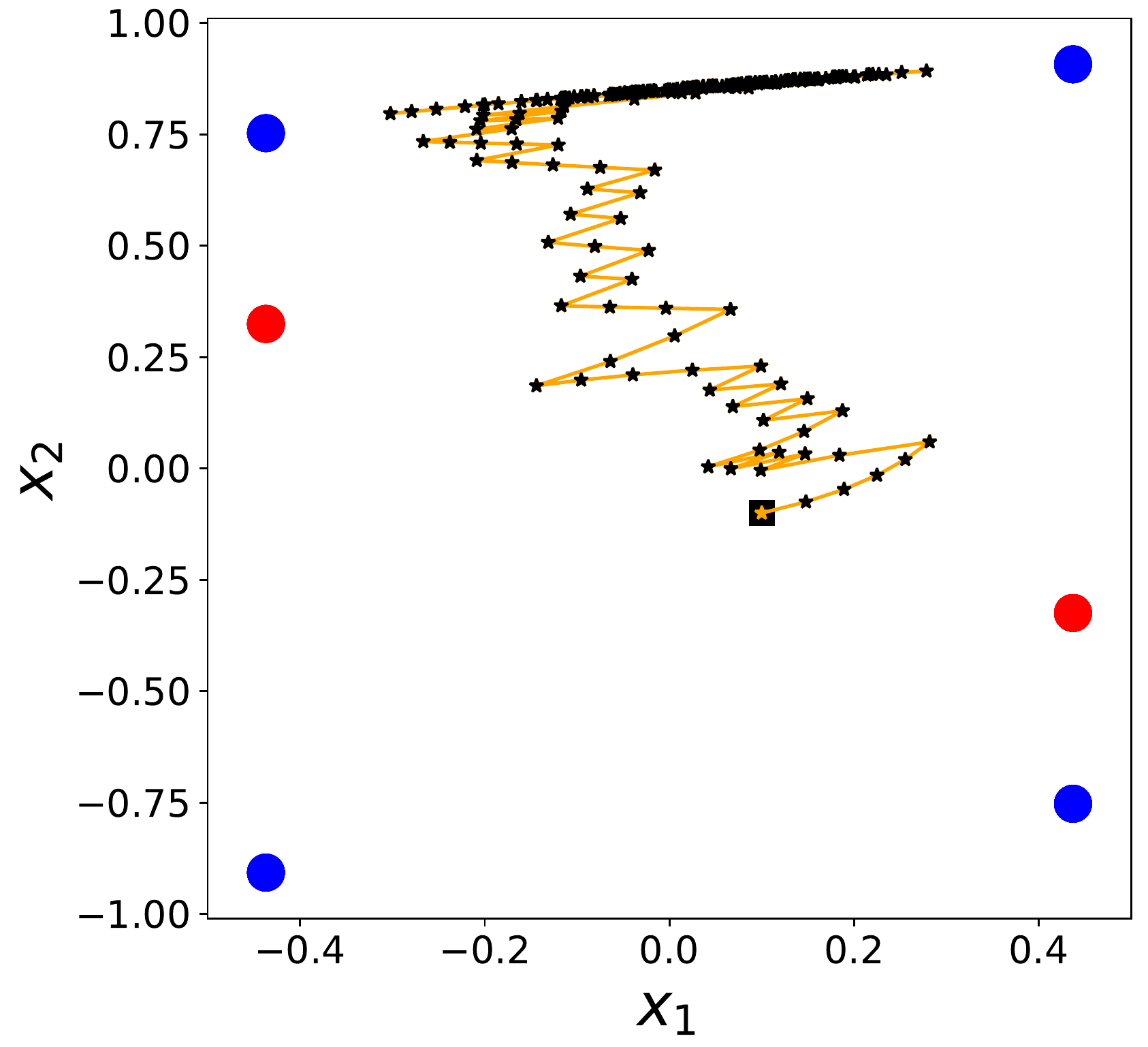}~\includegraphics[keepaspectratio=true,scale=0.45]{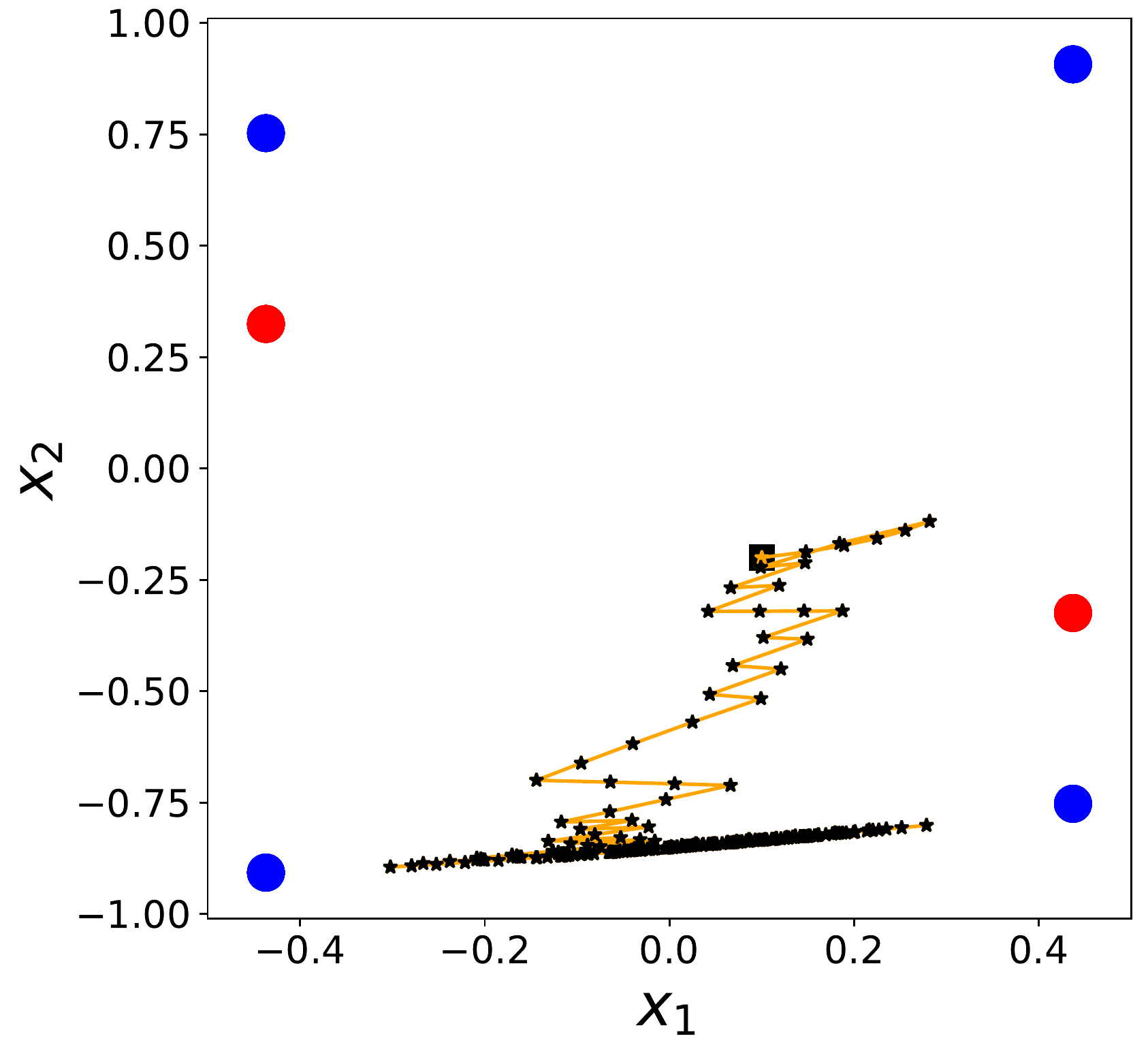}\\
\includegraphics[keepaspectratio=true,scale=0.31]{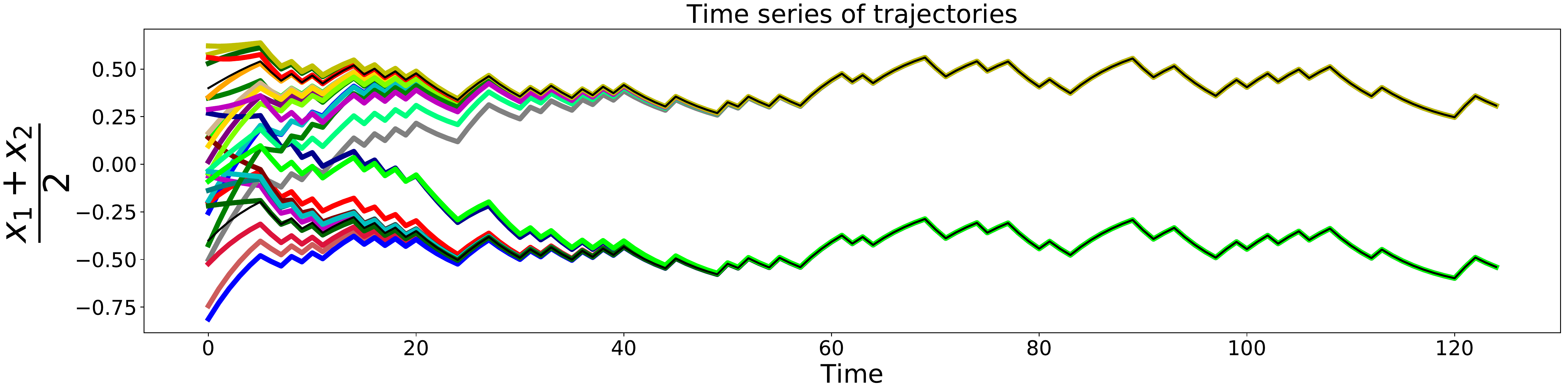}
\end{center}
\caption{
All pictures refer to the same input sequence $\bv \in \{ u_1, u_2 \}^{\bZ}$ which has been randomly generated with equal probability to switch between $u_1$ and $u_2$. 
\textbf{Top:} Fixed points of both (autonomous) component maps $f_1(x), f_2(x)$ are showed as blue (stable nodes) and red (saddle) dots, respectively. Notably, fixed points vertically aligned on the positive side of the $x_1$ variable characterise the $f_1$ map, while the ones on the negative side characterise $f_2$.
On the left panel, the initial condition of coordinates $(0.1,-0.1)$ evolves towards the upper UAES solution; on the right panel, the initial condition of coordinates $(0.1,-0.2)$ evolves towards the bottom UAES solution.
Under the same driving input sequence, almost every initial condition converges to one of those two nonautonomous attractors.
\textbf{Bottom:} Time series of the observable $ \frac{x_1+x_2}{2}$ of 30 trajectories starting from randomly chosen initial conditions are shown. 
This plot suggests a decomposition for the phase space in 2 UAESs, which leads the system to synchronise around two stable responses when driven by the given input sequence.
}
\label{fig:2asymp-solutions}
\end{figure}

The autonomous system $x[k+1] = f_1(x[k])$ has two asymptotically stable fixed points with a saddle between them along the vertical line of $ x_1 \approx 0.45 $, see Figure \ref{fig:2asymp-solutions}. Analogously, the autonomous system $x[k+1] = f_2(x[k])$ has two asymptotically stable points with a saddle between them along the vertical line $ x_1 \approx -0.45 $.

Exploiting Theorem \ref{thm:idea_contraction} we are able to prove the existence of two UAESs for every deterministic input sequence $\bu \in \{ u_1, u_2 \}^\bZ$. 
The Jacobian matrix of $f_1$ reads 
$$
    D_x f_1(x_1,x_2)= 
    \begin{bmatrix}
    1-\alpha/2\bigl[ 1 + \tanh^2(x_1/2 + 1/4) \bigl] & 0 \\
    0 & 1+\alpha/2\bigl[ 1 - 3\tanh^2(3x_2/2 + 3/20)\bigl]
    \end{bmatrix}.
$$
Diagonal elements are thus eigenvalues and also singular values, hence $\sigma( D_x f_1(x_1,x_2)) = 1+\alpha/2\bigl[ 1 - 3\tanh^2(3x_2/2 + 3/20) \bigl]$. 
Thus $ \lVert  D_x f_1(x_1,x_2) \rVert = 1+\alpha/2\bigl[ 1 - 3\tanh^2(3x_2/2 + 3/20) \bigl] $. 
The region of phase space where contraction occurs, i.e. $\lVert  D_x f_1(x_1,x_2) \rVert < 1 $, consists of 2 connected components divided by the strip $ -0.54 < x_2 < 0.34 $, approximately.
Similarly, for $f_2$ we have that $\lVert  D_x f_2(x_1,x_2) \rVert < 1 $ holds true everywhere except for $ -0.34 < x_2 < 0.54 $, approximately.
Now, it is easy to see that the region $R_+ = \{ (x_1, x_2) \in [-1,1]^2 : x_2 > 0.54 \} $ is positively invariant for both maps $f_1$ and $f_2$. 
Analogously, the region $R_- = \{ (x_1, x_2) \in [-1,1]^2 : x_2 < -0.54 \} $ is positively invariant for both maps $f_1$ and $f_2$. 
Hence, both regions $R_+$ and $R_-$ result to be convex $\cU$-positively invariant compact sets of the nonautonomous system where contraction occurs.
Therefore, Theorem \ref{thm:idea_contraction} ensures that inside the set $ R_+ $, there exists a UAES and, similarly, there exists another UAES inside $ R_- $.

For the simulation in Figure \ref{fig:2asymp-solutions}, we generated an input sequence with same probability to occur either $u_1$ or $u_2$.
In the top panel of Figure \ref{fig:2asymp-solutions}, we show the UAESs of such input-driven system. 
In the bottom panel, time series of the observable $\dfrac{x_1+x_2}{2} $ of many initial conditions are shown.
As can be seen, a UAES is substantially different from a fixed point due to its input-driven nature.

\clearpage
\subsubsection{The separatrix entire solutions}
\label{sec:separatrix}

\begin{figure}[ht!]
\begin{minipage}{\textwidth}
\centering
    \includegraphics[width=0.45\textwidth]{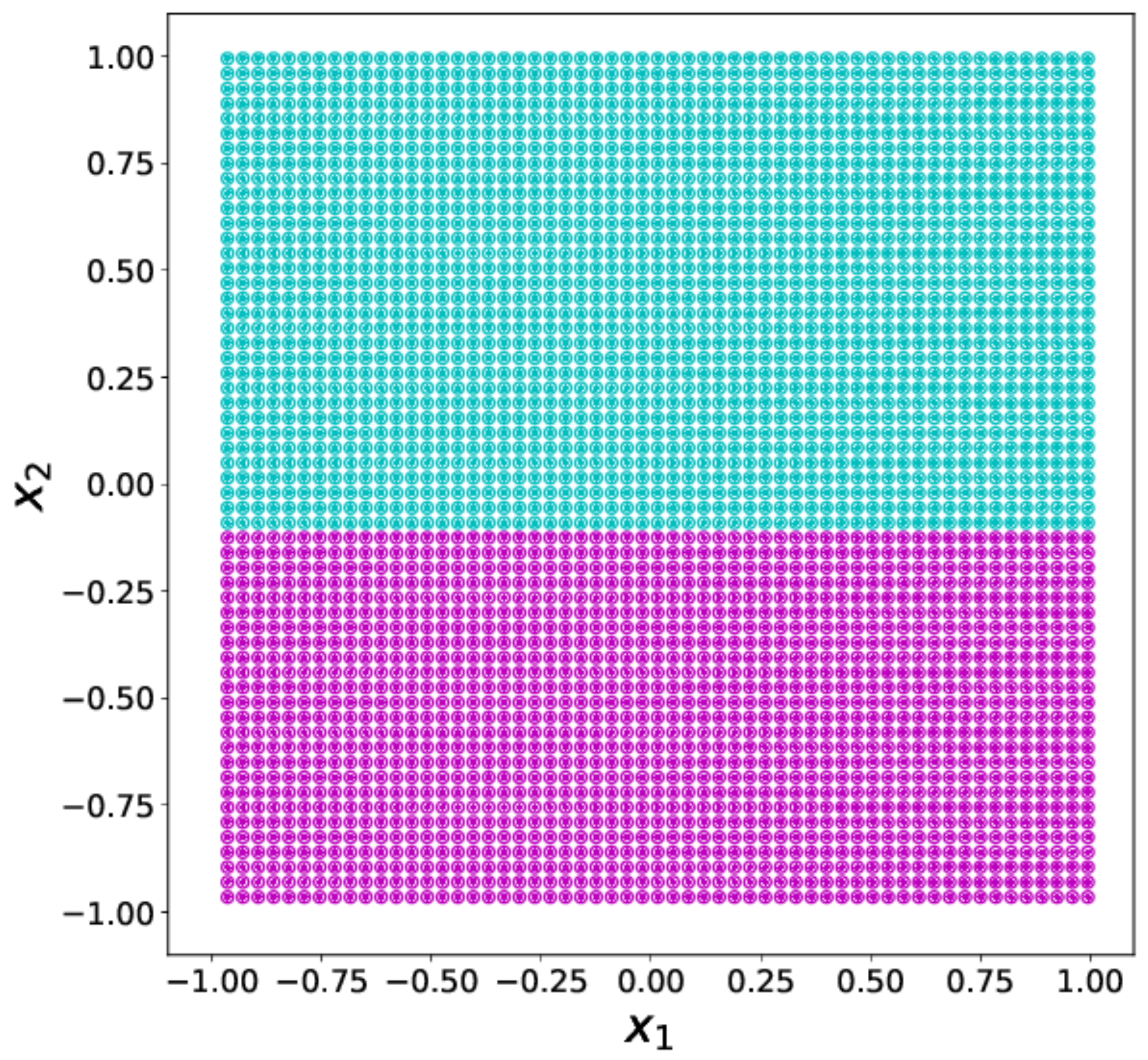}~\includegraphics[width=0.45\textwidth]{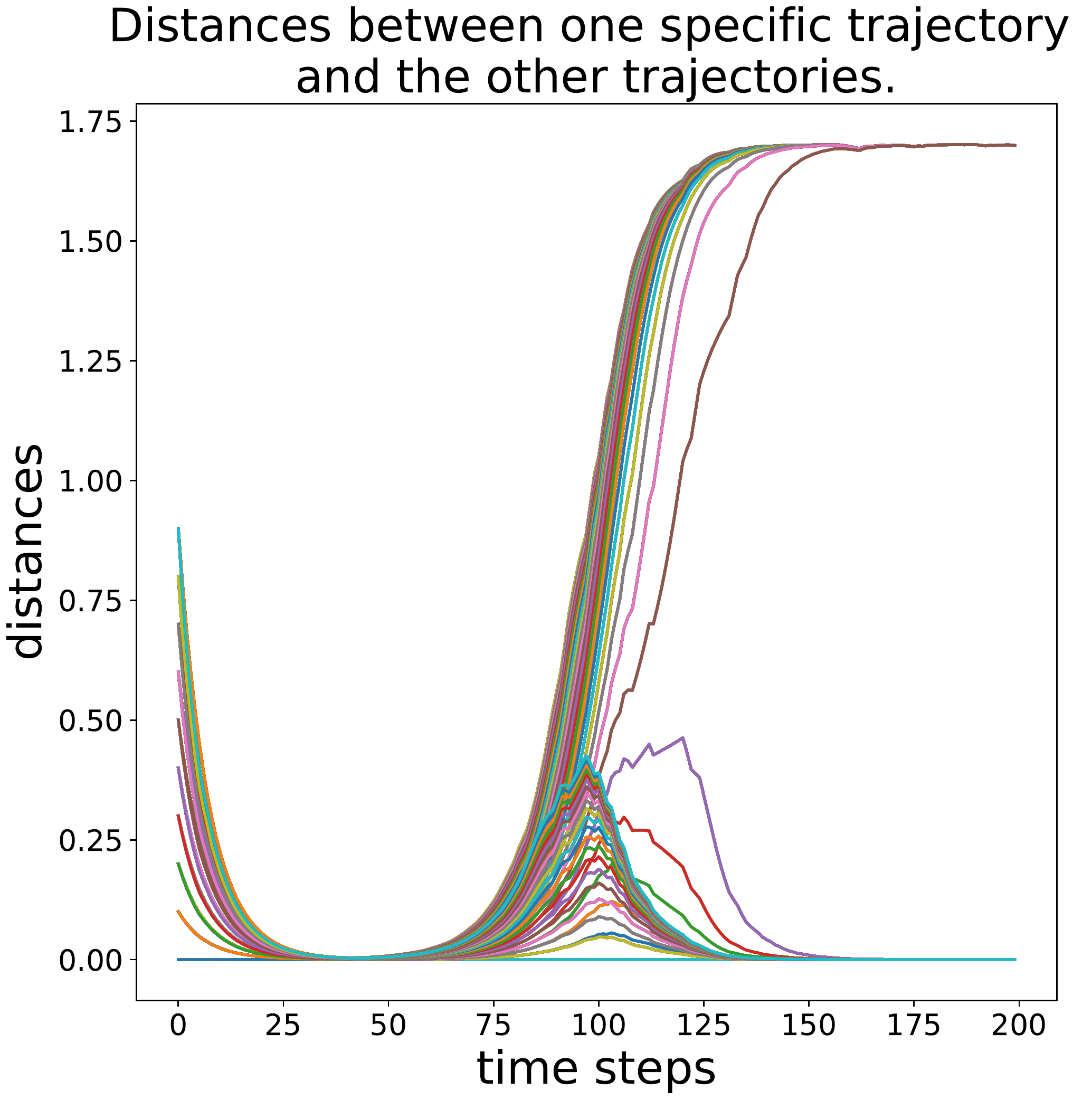}
\end{minipage}\\
\begin{minipage}{\textwidth}
\centering
    \includegraphics[width=0.45\textwidth]{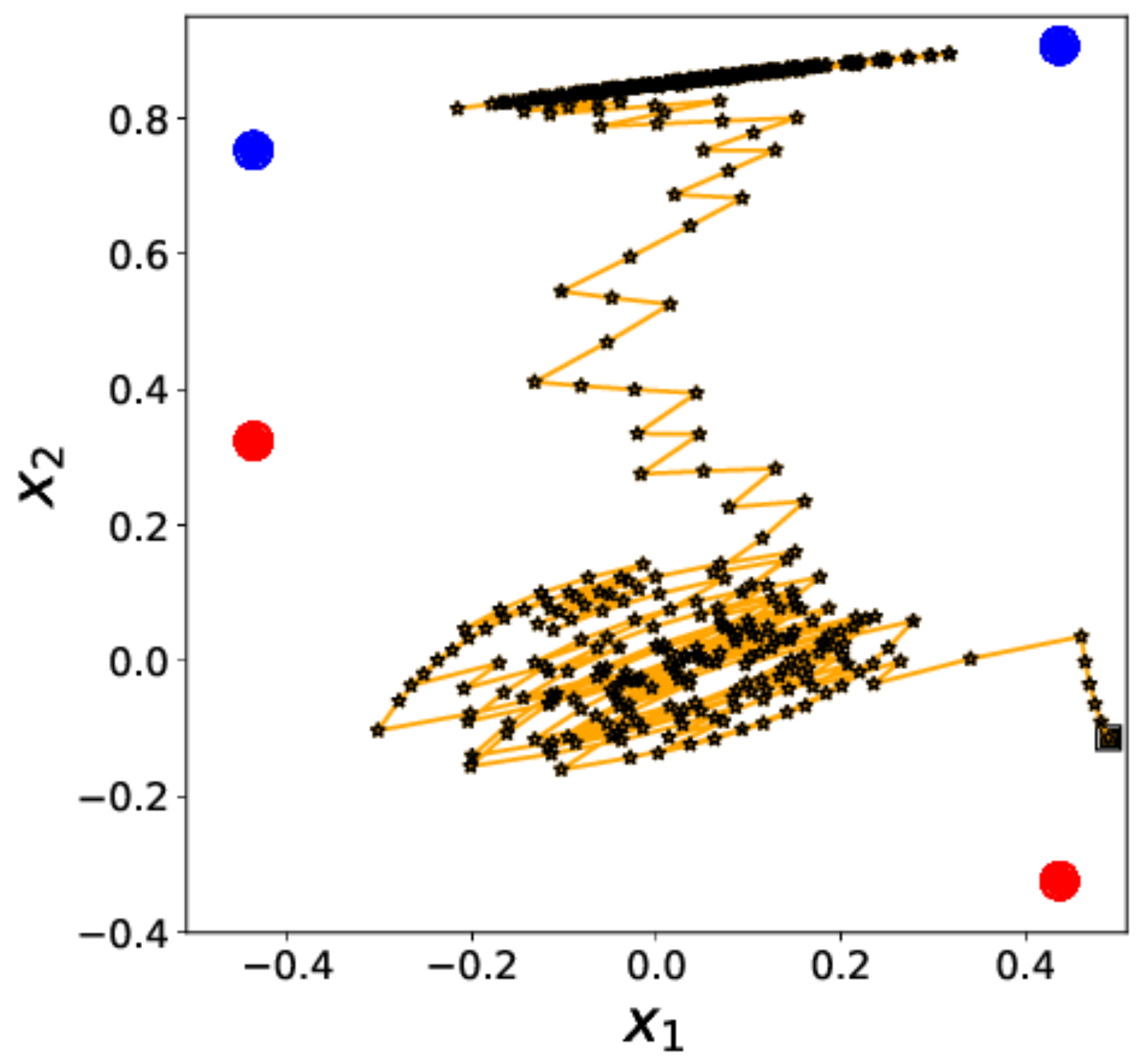}~\includegraphics[width=0.45\textwidth]{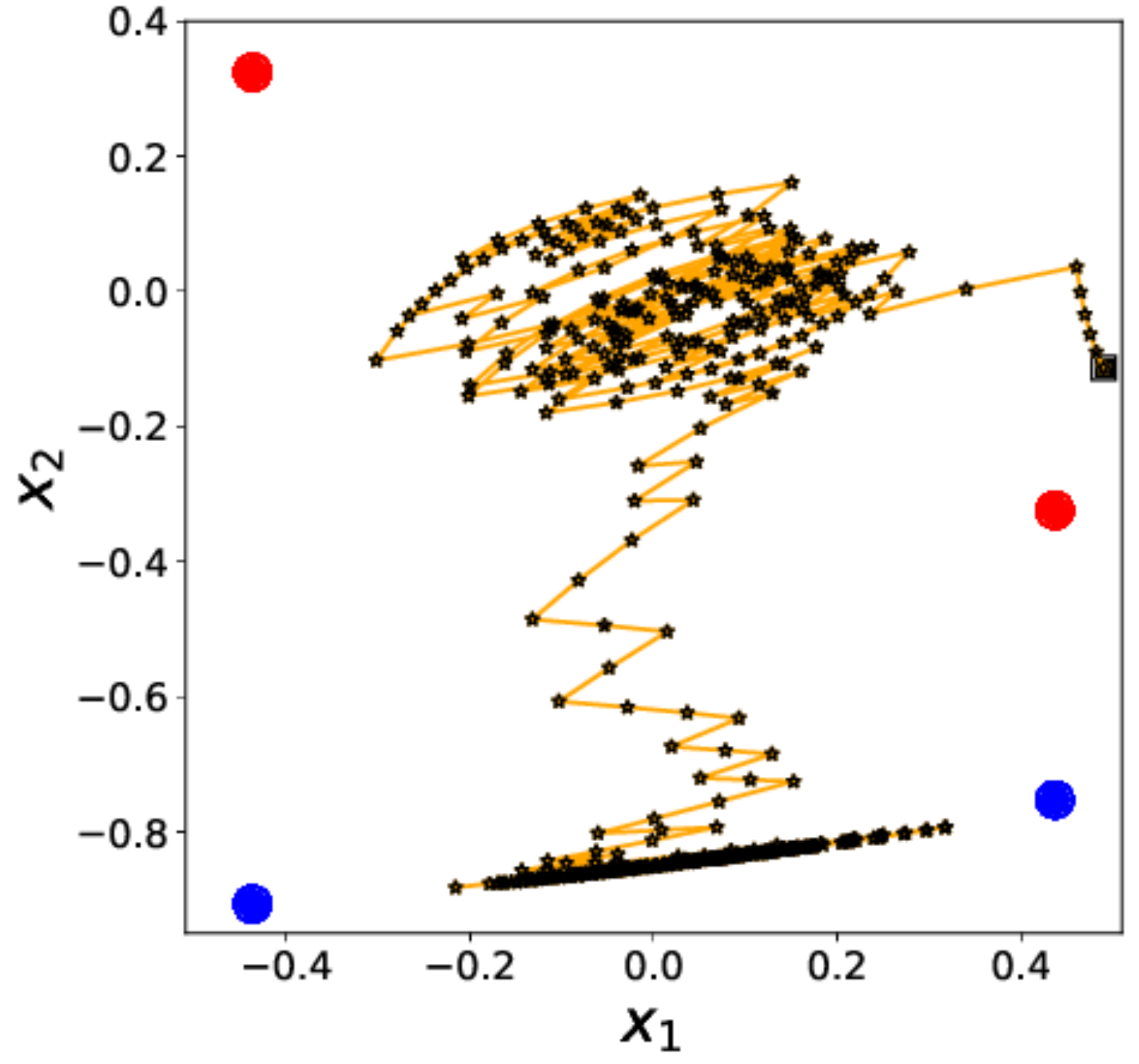}
\end{minipage}
\caption{
All pictures refer to the same randomly generated input sequence of Figure \ref{fig:2asymp-solutions}.
\textbf{Top left:} computation of the basins of attraction of the two UAESs at time step $k=0$.
\textbf{Top right:} 950 initial conditions have been chosen around the position of the separatrix at time step $k=0$; the $x_2$ coordinate starting from $-0.113248$ and cumulatively increasing of $10^{-6}$ until the value $-0.113198$, the $x_1$ coordinate starting from $-1$ and cumulatively increasing of $10^{-1}$ until the value $1$. The trajectory starting from the initial condition $(0, -0.113218)$ has been chosen as reference. At each time step, it has been computed and plotted the distance between the reference trajectory and all other trajectories.
\textbf{Bottom left:} initial condition is $(0.49, -0.11322366651)$. The trajectory rambles in phase space following the separatrix attractive entire solution, then after about 230 time steps it starts to converge towards the upper UAES.
\textbf{Bottom right:} here the initial condition is $(0.49, -0.11322366652)$, hence it differs of $10^{-11}$ from the one in the left picture. The trajectory practically coincides with the one in the left picture but eventually converges towards the bottom UAES.} 
 
\label{fig:separatrix} 
\end{figure}

The basin boundary of the two UAESs is, at each fixed time step, a horizontal line which lies in-between the stable manifolds of the two saddles of the component maps.
In the top left of Figure \ref{fig:separatrix}, this line is shown for time step $k=0$ for the particular input sequence used in simulation of Figure \ref{fig:2asymp-solutions}.
In the top right plot of Figure \ref{fig:separatrix}, we evolved many different initial conditions in a neighbourhood of such a horizontal line. We randomly chose one of these trajectories as reference trajectory and computed the distances, at each time step, between each trajectory with respect to a randomly chosen one.
It is numerically observed that the closer we start to such horizontal line, the more time is needed to converge towards one of the UAESs.
Indeed, the uniform convergence (i.e. in the Hausdorff semi-metric sense) \cite{kloeden2011nonautonomous,falconer2004fractal} towards a UAES holds only excluding a neighbourhood of the boundary of its basin of attraction.
Moreover, from the top right plot of Figure \ref{fig:separatrix} we observe that each initial condition seems to be attracted for a few time steps to a special trajectory (all the distances go to zero), but are eventually pushed away from it in the long run towards one of the two UAESs.
This suggests the existence of a third attractive entire solution whose basin of attraction is exactly the (moving in time) separatrix line.
As a matter of fact, each point of such separatrix line represents a fibre of an entire solution. Interestingly, all those entire solutions seem to converge towards a unique solution which appear as in bottom pictures of Figure \ref{fig:separatrix}. Note that, in particular, such nonautonomous system admits an infinite number of entire solutions but only two of them essentially characterise its dynamics\footnote{Note that all the infinitely many separatrix entire solutions are part of the global pullback attractor of Definition \ref{def:pullback_attractor}. Moreover, all the entire solutions forward converging to the two UAESs and originating in the infinite past from the separatrix solutions are also them part of the global pullback attractor of Definition \ref{def:pullback_attractor}.}.
Although this third attractive entire solution plays a key role in the nonautonomous dynamics of the system, it should not be considered equivalent to the other two UAESs. In fact, the probability to pick up an initial condition in phase space converging to this particular entire solution is null. 
This is the reason why we require for a UAES to have a neighbourhood of attraction, as formally stated in Definition \ref{def:local_point_attractor}, which in particular implies that its basin of attraction has positive Lesbegue measure.
In the bottom panels of Figure \ref{fig:separatrix}, we tried to detect the third attractive entire solution letting evolve initial conditions extremely close to the separatrix line. Such a special entire solution appears unstable and it traces an orbit that wanders erratically in the region between the stable manifolds of the two saddles of $f_1$ and $f_2$.

\clearpage
\subsection{An example with input-dependent echo index}
\label{sec:Change-echo-index}

In this section, we consider a one-dimensional RNN and show how the echo index may vary according to the specific input sequence driving the dynamics.

Let us consider the input-driven system $x[k+1]=G(u[k+1],x[k])$, with 
\begin{equation}
    \label{eq:oneDim_rnn}
    G(u,x) = \tanh( 1.01 x + w u ), \qquad x \in [-1,1], \quad u \in [-1,1],
\end{equation}
where $w$ is a positive real value playing the input gain role.
Note that in this example the value $1.01$ represents the maximum singular value (and spectral radius as well) of a one-dimensional recurrent layer.
Therefore, the autonomous system defined by the map $F(x):=G(0,x)$ is expanding around the (unstable) fixed point $x=0$ and there exist two (uniformly attracting) stable points $x_1^* \approx -0.17 , x_2^* \approx 0.17$.

Theorem \ref{thm:small_input} and Proposition \ref{prop:large_inputs} suggest that the echo index is also determined by the amplitude of the inputs driving the dynamics. Here, we show how modulating the input amplitude via $w$ yields different values of echo index.
Note that, as $w$ scales the inputs provided to the system \eqref{eq:oneDim_rnn}, in practice we force the autonomous map $F$ with inputs assuming values in $[-w,w]$.
It is possible to analytically compute the value of $w$ such that a fold bifurcation occurs in the system \eqref{eq:oneDim_rnn}, which is approximately $ w_c \approx 0.0007 $; see Eq. 36 of \cite{ceni-ff-18} for details. This implies that perturbing the autonomous dynamics by means of any input sequence with amplitude less than $w_c$ will induce an input-driven system with echo index 2.

\begin{figure}[ht!]
\begin{center}
\includegraphics[keepaspectratio=true,scale=0.37]{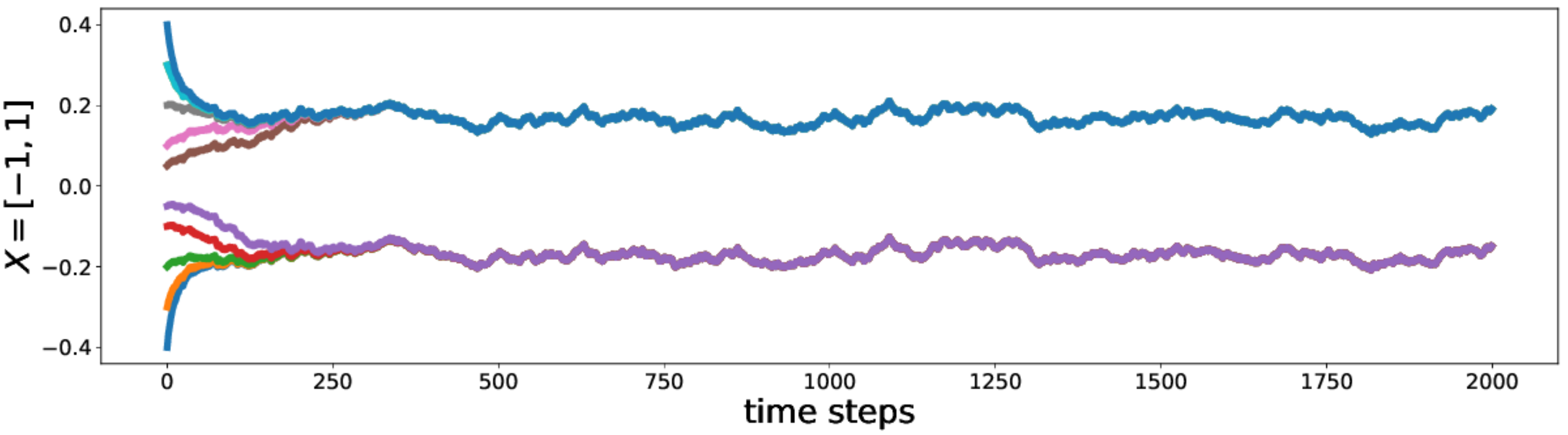}\\
\includegraphics[keepaspectratio=true,scale=0.37]{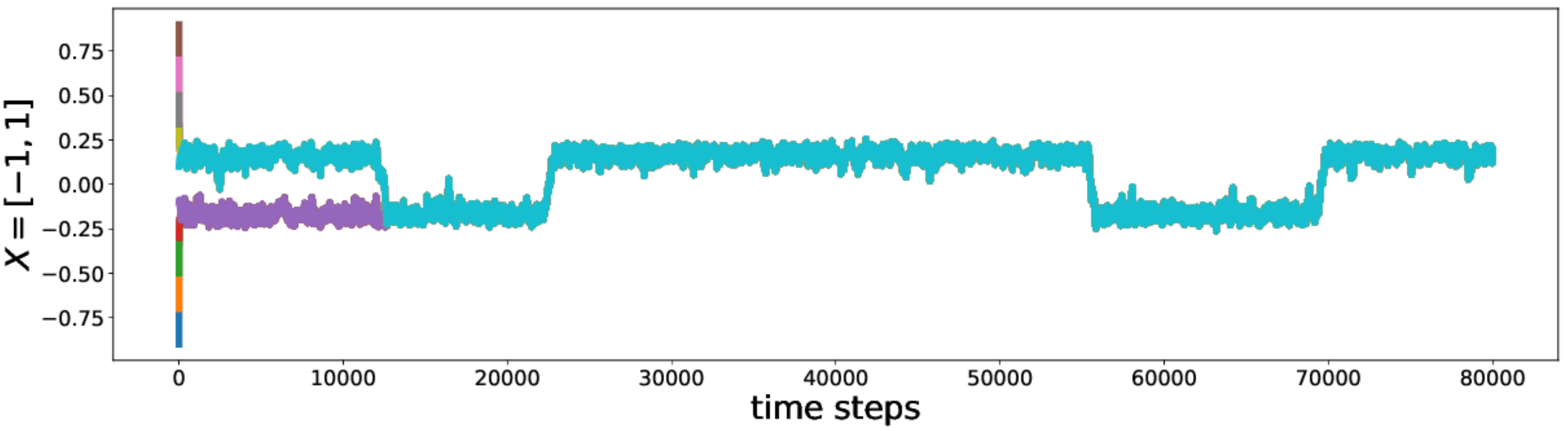}\\
\includegraphics[keepaspectratio=true,scale=0.37]{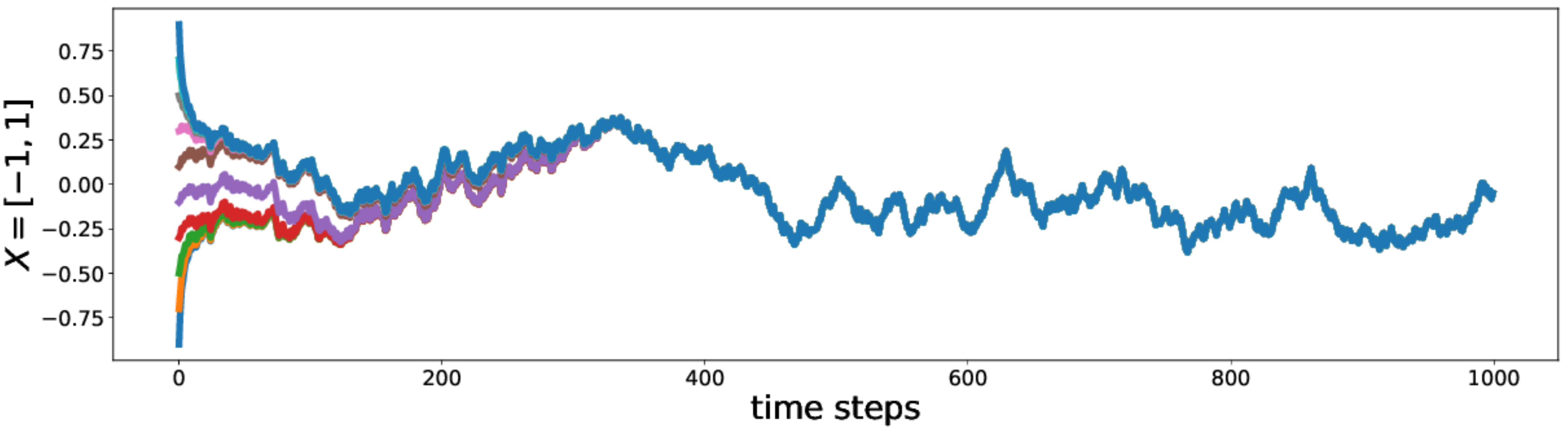}
\end{center}
\caption{
Nonautonomous dynamics of the system \eqref{eq:oneDim_rnn} driven by an input sequence generated according to a uniform distribution in $[-1,1]$ and then scaled by means of a positive real parameter $w$, for three values of the parameter $w = 0.0006, 0.01, 0.05$.
Ten initial conditions have been run in all three cases, represented with different colours in the plots.
\textbf{Top:} case of $w=0.0006$. The input driven system exhibits echo index 2, i.e. there are two UAESs.
\textbf{Centre:} case of $w=0.01$. The input driven system has echo index 1. There exists a unique UAES whose behaviour is affected by the vicinity of the fold bifurcation of the underlying autonomous map $F(x)=G(0,x)$. The resulting dynamics manifest a switching motion. 
\textbf{Bottom:} case of $w=0.05$. The input-driven system presents echo index 1 and the corresponding uniformly attracting entire solution does not exhibit any switching behaviour.
}
\label{fig:1Dmap}
\end{figure}

To show evidence of this, we generated an input sequence according to a uniform distribution in $[-1,1]$, then we scaled it by means of three values of $w = 0.0006, 0.01, 0.05$.
Figure \ref{fig:1Dmap} shows the evolution of ten initial conditions of the system for the three values of $w$.
In the first case of $w=0.0006$ (top Figure \ref{fig:1Dmap}), the input-driven system exhibited echo index 2, in accordance with the fact that the critical value $w_c$ is greater than $w=0.0006$.
The second case of $w=0.01$ (centre Figure \ref{fig:1Dmap}) produced an interesting dynamics with echo index 1: the unique UAES is characterised by a switching behaviour between the ``ghosts'' of the two randomly perturbed autonomous stable solutions of $F$.
Finally, in the third case of $w=0.05$ (bottom Figure \ref{fig:1Dmap}) the nonautonomous dynamics converge toward a unique UAES (i.e. echo index 1), which wanders randomly around the origin.

\clearpage
\subsection{RNNs dynamics in a context-dependent task}
\label{sec:Maass}

Here, we train an RNN to solve the task described in \cite[Figure 5]{hoerzer2012emergence}, which consists of performing some context-dependent computation.
Our results suggest that the RNN dynamics are characterised by a decomposition in two UAESs, i.e. the trained RNN has echo index 2.

The task is the following.
The RNN dynamics is driven by a bi-dimensional input sequence, $ \bu = \{ (u_1[k], u_2[k]) \} $, and produces bi-dimensional output $ \mathbf{z}= \{ (z_1[k], z_2[k]) \}$. The behaviour of the input-driven RNN dynamics is controlled by two impulsive control inputs, $u_3[k], u_4[k]$, which do not contribute in driving the dynamics but give the RNN a context of the ``on'' and ``off'' type.
These instantaneous pulses have unitary amplitude and occur with probability $0.01$, i.e. $u_3[k], u_4[k]$ are null most of the time and occur, on average, every 100 time steps.
The first readout is used to produce an output $z_1[k]$ which maintains a working memory of the context: assuming the value $+1$ for the ``on'' state and $-1$ for the ``off'' state.
For this reason, the first output $z_1[k]$ is fed back into the network via the output feedback connections, $W_{fb}$.
The two input sequences driving the dynamics, $u_1[k], u_2[k]$, are two independently generated time-varying signals obtained through the discrete convolution\footnote{Given two sequences $ \{ a[k] \}_{k \in \bZ}$ and $ \{ b[k] \}_{k \in \bZ}$, the discrete convolution is defined as $(a * b)[n] = \sum_{m = -\infty}^{\infty} a[m] b[n - m]$.} of a uniformly distributed signal assuming values in $[0,1)$ with the smooth exponential filter $ g(s) = \exp(-s/50 ) $.
Bias values of, respectively, $0.3$ and $0.15$ were added to these convolutions; input sequences $ u_1[k], u_2[k]$ are normalised so their maximum value is one.
The second readout has to learn the routing of $u_1[k]$ to output $z_2[k]$ if the network is in state ``on'', and $u_2[k]$ to output $z_2[k]$ if the network is in state ``off''.

In Figure \ref{fig:Maass} we report the results of a trained RNN with $N_r=500$ neurons; in Appendix \ref{sec:training}, we provide details regarding how we trained a RNN of the form \eqref{eq:leaky_rnn}-\eqref{eq:output} to solve this task.
In the top two plots of Figure \ref{fig:Maass}, the outputs $z_1[k], z_2[k]$ produced by the RNN versus the target signals are shown, demonstrating that the network learned to solve the task quite well.
The three plots in the middle of Figure \ref{fig:Maass} show the evolution of the test trajectory projected in the two-dimensional subspace of the RNN phase space obtained by computing the first two principal components of the state trajectory.
In the two plots at the bottom of Figure \ref{fig:Maass}, we switch off the control inputs and run 100 initial conditions uniformly distributed in phase space. These last two plots show the presence of two UAESs describing the behaviour of the RNN trained to solve such a context-dependent task.
\begin{figure}[]
\begin{minipage}{\textwidth}
    \centering
    \includegraphics[keepaspectratio=true,scale=0.25]{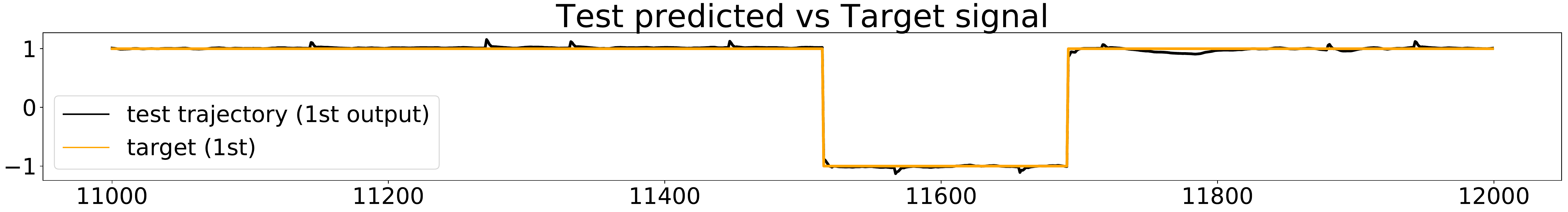}\\
    \includegraphics[keepaspectratio=true,scale=0.25]{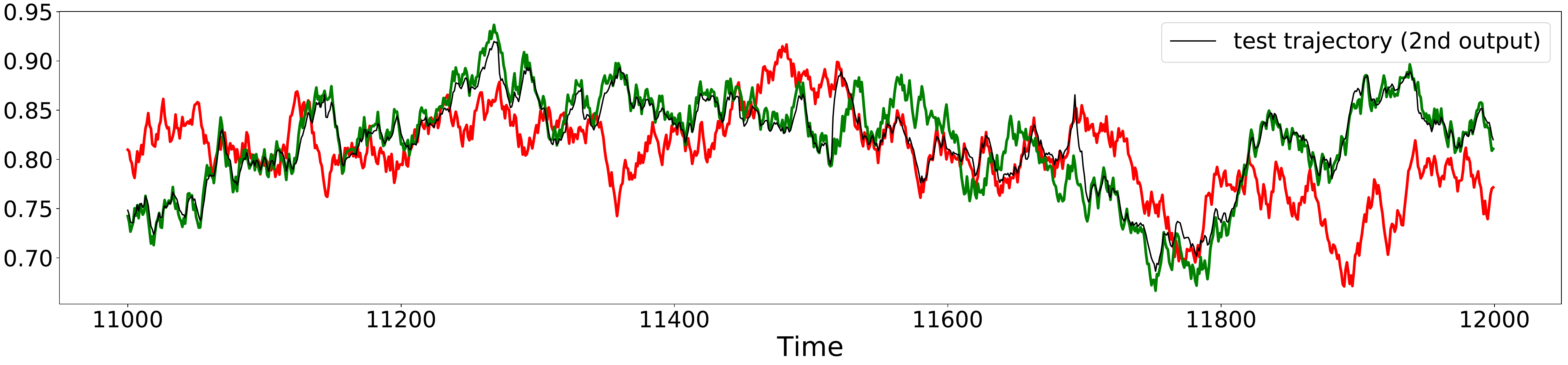}\\
    \includegraphics[keepaspectratio=true,scale=0.33]{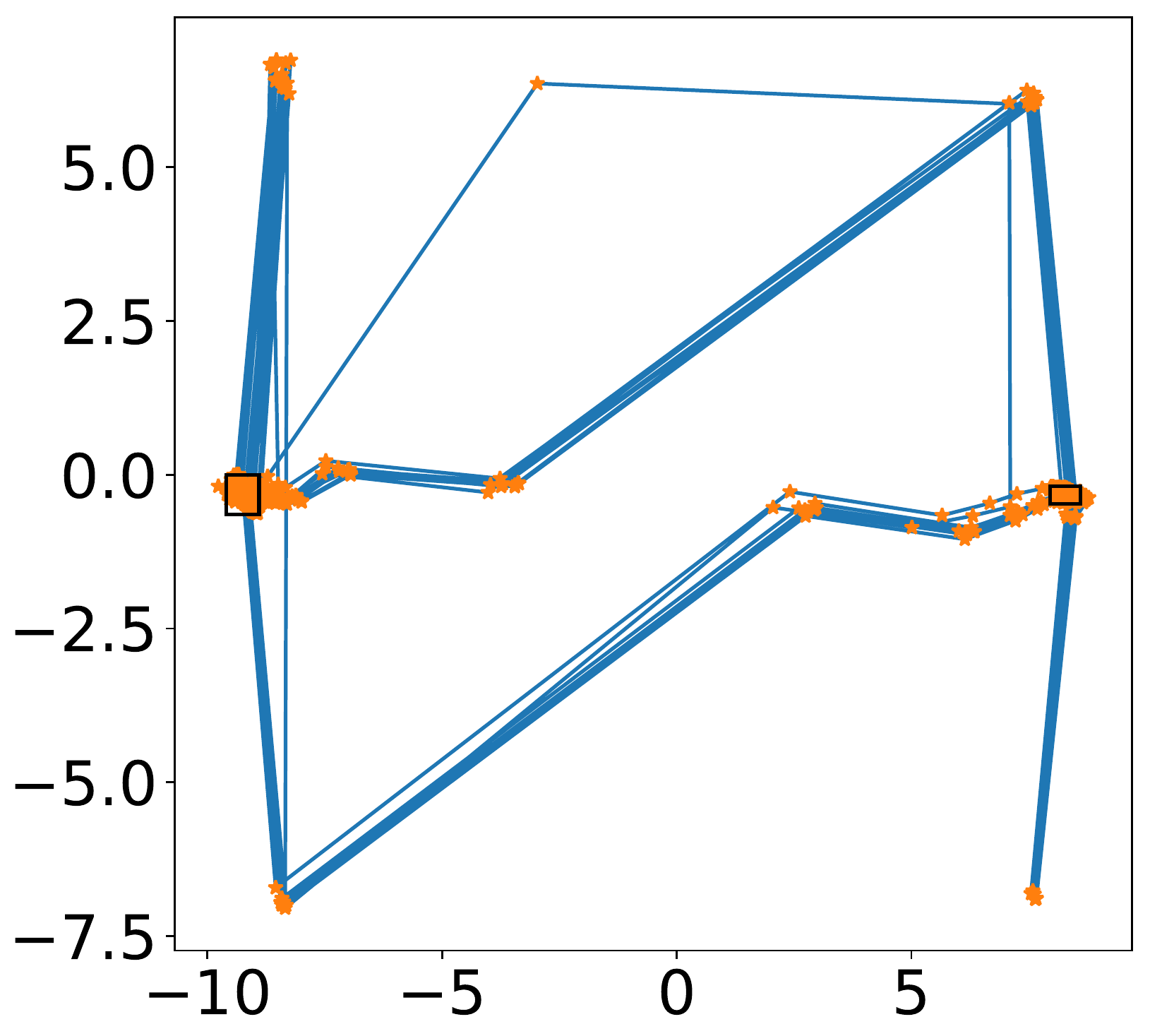}~\includegraphics[keepaspectratio=true,scale=0.33]{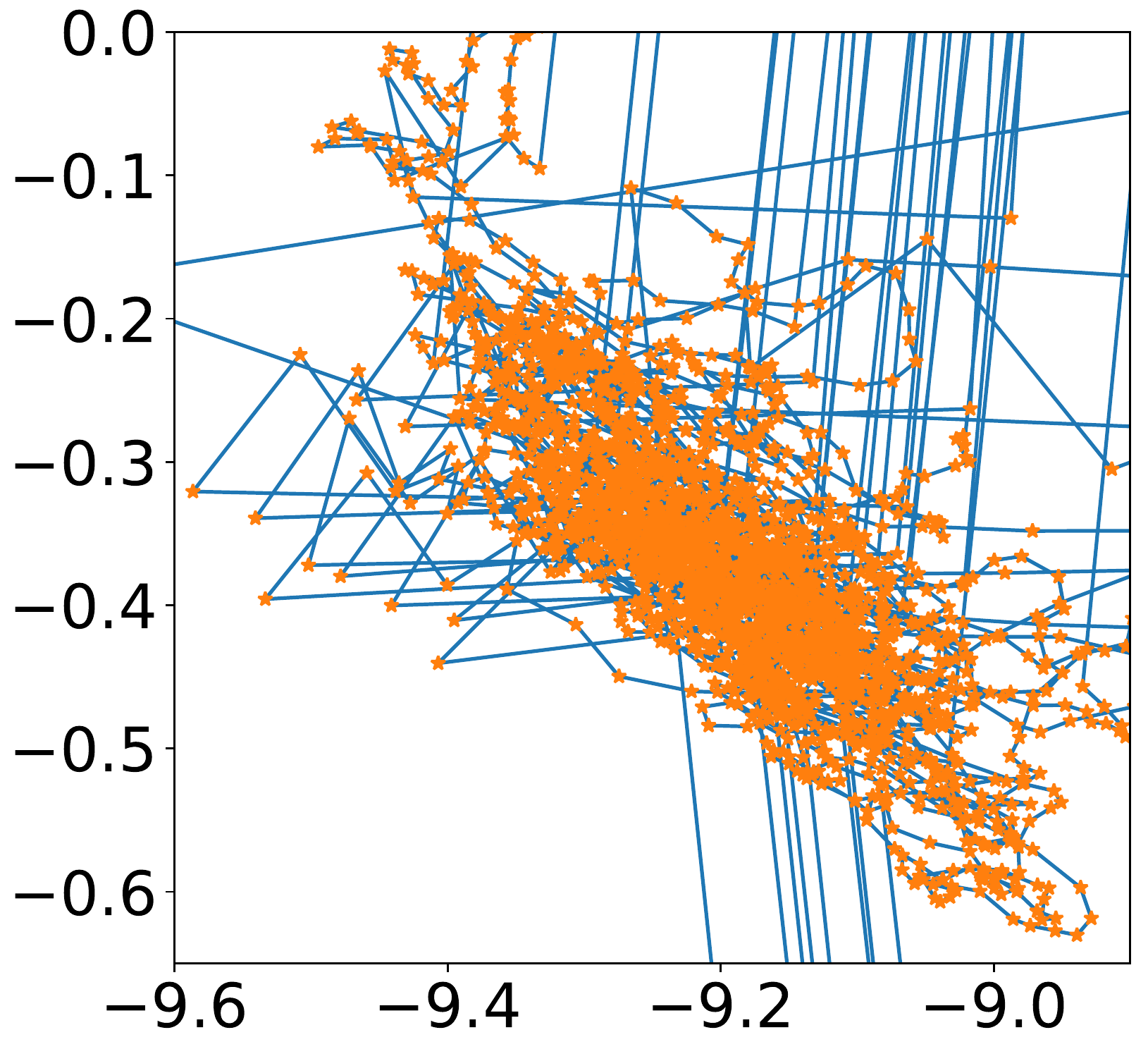}~\includegraphics[keepaspectratio=true,scale=0.33]{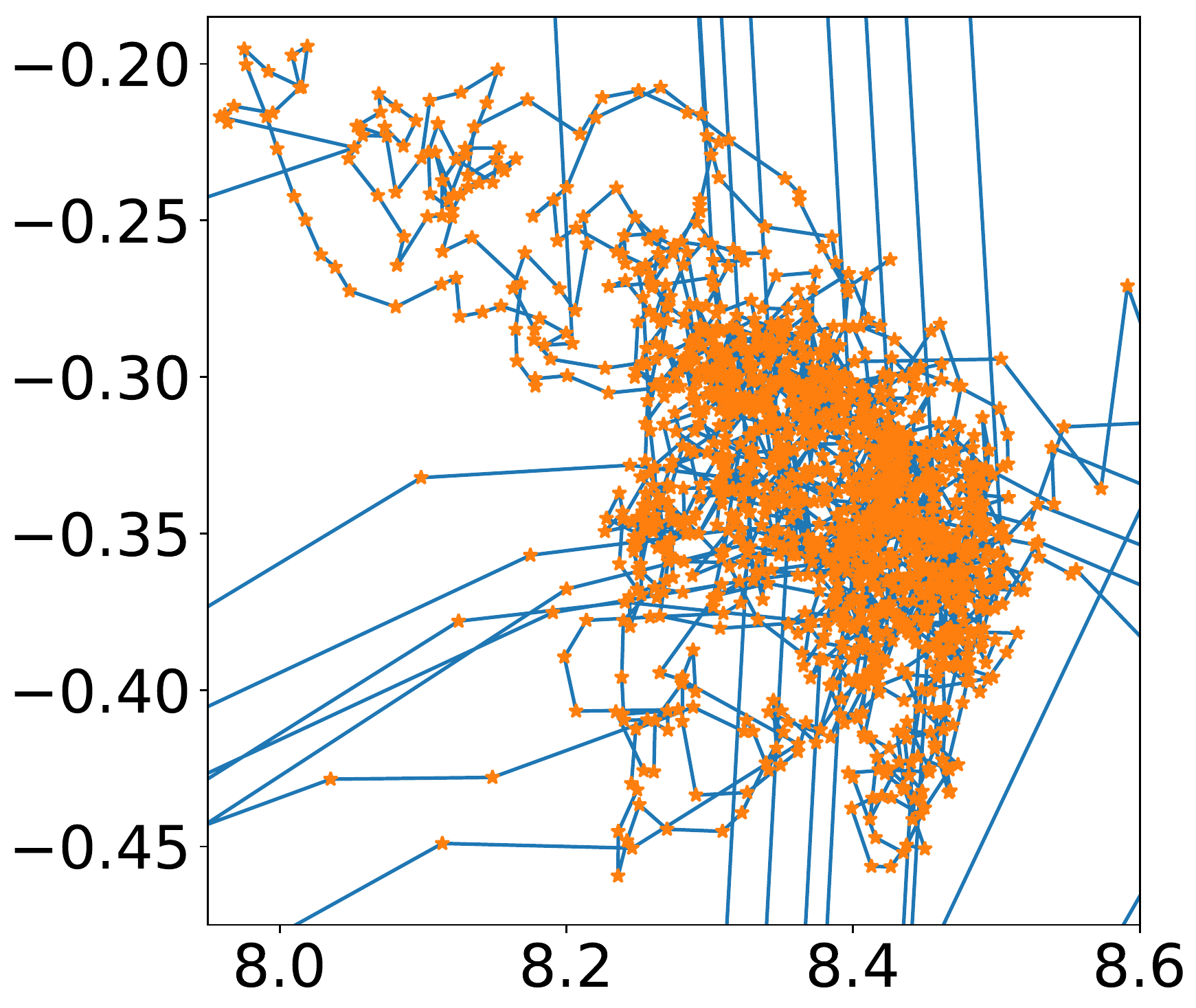} \\
    \includegraphics[keepaspectratio=true,scale=0.25]{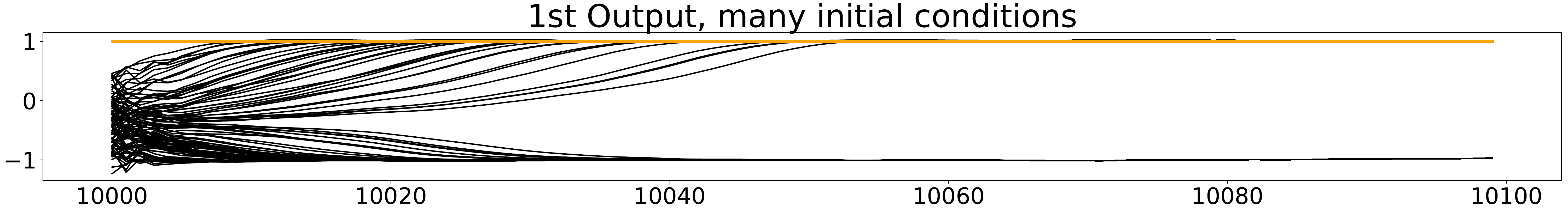}\\
    \includegraphics[keepaspectratio=true,scale=0.25]{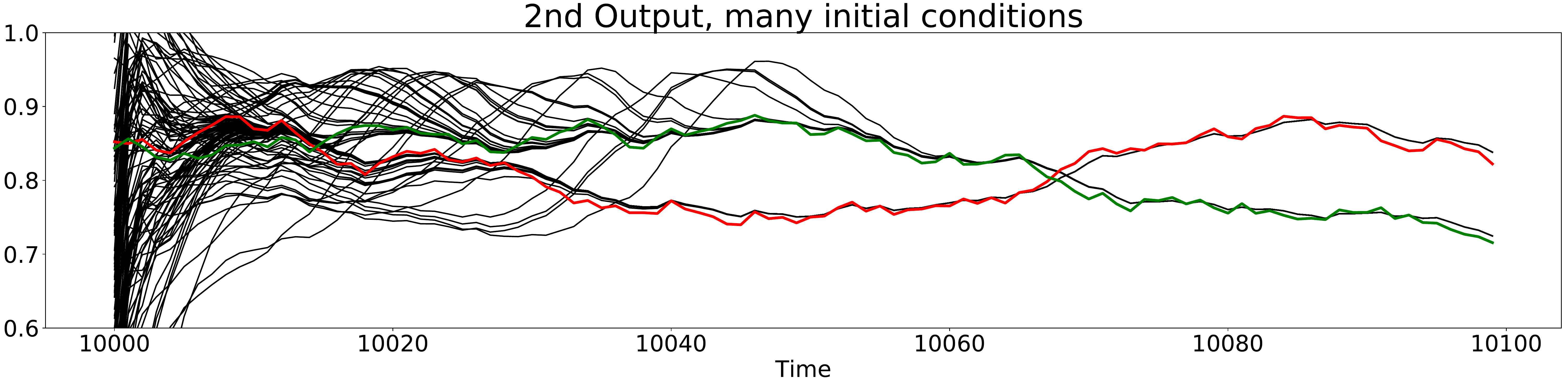}\\
\end{minipage}
\caption{
All plots refer to the test session, after training is completed.
The internal context-state to be learned is shown in orange and assumes values $+1$ for the ``on'' state or $-1$ for the ``off'' state, while the inputs $u_1[k], u_2[k]$, driving the dynamics are shown in green and red, respectively. 
\textbf{Top:} the top two plots show the outputs $z_1[k], z_2[k]$ (in black) versus the target for the first output (in orange) and for the second output (in green when the state is ``on'', and in red when the state is ``off'').
\textbf{Middle:} the left middle plot depicts the evolution of the test trajectory projected in the two-dimensional space spanned by the two principal components of the network states. The cumulative variance of the first two principal components is $0.98$.
The test trajectory spends most of the time around two locations delimited in the figure by black frames.
The other two plots in the middle show a zoomed view of the test trajectory inside of these black frames.
The attractor on the left corresponds to the ``on'' state and it is related to the green signal, while the attractor on the right corresponds to the ``off'' state and it is related to the red signal.
While these two nonautonomous attractors look very complex and unstable, they are in fact very stable to variations of initial conditions, as demonstrated in the bottom plots.
\textbf{Bottom:} Evolution of the trained RNN dynamics with 100 initial conditions randomly distributed in phase space. All initial conditions converge either to the left or to the right nonautonomous attractor.
This experiment strongly suggests that the RNN training produced a phase space decomposition in two UAESs, which in turn are used to solve the context-dependent task.
}
\label{fig:Maass}
\end{figure}

\clearpage
\section{Conclusions} 
\label{sec:conclusions}

As main contribution of this paper, we highlight how the echo state property, which guarantees the existence of a unique (stable) response to an input sequence, may be generalised so that only local behaviour in phase space in taken into account. Accordingly, we show how a recurrent neural network might reliably produce several stable responses to an input sequence: the echo index introduced here counts such stable responses.
Our theoretical developments are framed within the theory of nonautonomous dynamical systems. We introduce a suitable definition of attractor for input-driven recurrent neural networks that we call a uniformly attracting entire solution, which models and characterises responses of recurrent neural networks to input sequences.
The presence of more than one stable response indicates the possibility to observe and exploit multiple, yet consistent behaviours of a recurrent neural network driven by an input sequence (for example, the context-dependent task shown in Section~\ref{sec:Maass} requires learning two stable responses).
On the other hand, echo index greater than one might also indicate incorrect training and hence signal possible malfunctionings on a task requiring a unique behaviour in phase space.

We believe that the notions and results introduced here will prove fundamental for the more general and ambitious goal of providing mechanistic models describing the behaviour of recurrent neural networks in machine learning tasks, such as time series classification and forecasting.
Moreover, our developments highlight the need to rectify what is commonly misunderstood as unreliable behaviour of driven RNNs with the occurrence of multistable, nonautonomous dynamics. This paradigm shift will ultimately lead to a more suitable definition of ``chaotic behaviour'' for driven RNNs and related onset mechanisms.

\section*{Acknowledgements}

LL gratefully acknowledges partial support of the Canada Research Chairs program. CP and PA gratefully acknowledge partial support from the NZ Marsden fund, grant UOA1722. PA acknowledges funding from EPSRC as part of the Centre for Predictive Modelling in Healthcare grant EP/N014391/1.
The authors thank the anonymous reviewers for their helpful comments, and Manjunath Gandi for helpful discussions.


\bibliographystyle{elsarticle-num}

\bibliography{nesp-refs}

\clearpage
\appendices

\section{Hausdorff distance}
\label{sec:hausdorff}

Let $(X,d_{X})$ be a metric space: for convenience we recall some basic properties of the Hausdorff distance between two subsets of $X$.
For any subset $ Y \subset X $, we denote $ B_\varepsilon(Y) $ the $\varepsilon$-neighbourhood of the set $Y$, i.e. $B_\varepsilon(Y) := \bigcup_{y \in Y} B_\varepsilon(y) $ where $ B_\varepsilon(y) := \{ x \in X \,\, | \,\, d_{X}(x,y) < \varepsilon \} $ is the open ball of radius $\varepsilon$ centred on $ y \in Y$.
\begin{defn}
For any pair of nonempty subsets $ Y,Z \subset X $, let us define the following function 
\begin{equation}
\label{eq:semihausdorff}
    h(Y,Z) := \sup_{y \in Y}\underbrace{ \inf_{z \in Z}{d_{X}(y,z)}. }_{=:d_{X}(y,Z) \text{ point-set distance}}
\end{equation}
We call $ h :(\powerset(X)\setminus \{ \emptyset \})^2 \rightarrow [0, +\infty] $ the \emph{Hausdorff semi-distance} \cite{falconer2004fractal} of the metric space $(X,d_{X})$, where $\powerset(X)$ is the power set of $X$.
\end{defn}
Function $ h $ is not symmetric, i.e. in general $ h(Y,Z) \neq  h(Z,Y)$. Unfortunately, the fact that $ h(Y,Z)=0$ does not imply that $ Y = Z $. Nevertheless, it holds that $h(Y,Z)= \inf \{ \varepsilon \geq 0 \,\, | \,\,   Y \subseteq B_\varepsilon(Z)  \} $ hence the following is true for all $\varepsilon \geq 0 $,
\begin{equation}
    \label{eq:hausdorff-property}
    h(Y,Z)=\varepsilon \quad  \Longrightarrow \quad  Y \subseteq \overline{B_\varepsilon(Z)},
\end{equation}
where $\overline{B}$ denotes the closure of a set $B$.
Moreover, if the subset $Y$ is bounded then $ h(Y,Z) < +\infty $.
Therefore, the function 
\begin{equation}
    \label{eq:hausdorff}
    H(A,B) := \max \{ h(Y,Z) \, , \, h(Z,Y)  \} 
\end{equation}
is a metric on the space of all nonempty \emph{compact} subsets of $X$ and we call it the \emph{Hausdorff distance} of the metric space $(X,d_{X})$.

\section{Pullback attractors of input-driven RNNs}
\label{sec:pullback_thms}

Here, the state space is assumed to be a complete metric space and is denoted as $(Y,d_Y)$, with $ Y \subseteq \mathbb{R}^{N_r}$ closed and $d_Y$ is the Euclidean distance. Hence, we relax the hypothesis of compactness and deal with a state space that is not necessarily bounded.

For the purpose of this section, it will be useful to introduce the definition of pullback absorbing set.
In this regard, we wish to recall from Section \ref{sec:cocycle_formalism} that $\sigma : \cU  \longrightarrow \cU $ (the shift operator) defines an autonomous dynamical system acting on a compact metric space $(\cU,d_{\cU}) $ (the space of all admissible input sequences), and $\Phi: \bZ_0^+ \times \cU \times Y \longrightarrow Y $ (the cocycle map) describes the nonautonomous dynamics on a complete metric space $(Y,d_{Y})$.
\begin{defn}\cite[Definition 3.17]{kloeden2011nonautonomous}
\label{def:pullback_absorbing}
A nonempty compact subset $B \subseteq Y $ is called \emph{pullback absorbing} for a family of inputs $\cV\subseteq \cU $ if
$$
\forall \, \bu \in \cV, \, \forall \, \text{bounded } D\subseteq Y , \, \exists \, N=N(\bu,D)\in \bZ_0^+ \, : \, \Phi(n,\sigma^{-n}(\bu),D) \subseteq B \,\, \forall n\geq N.
$$
Analogously, a nonempty compact subset $B \subseteq Y $ is called \emph{forward absorbing} for a family of inputs $\cV\subseteq \cU $ if
$$
\forall \, \bu \in \cV, \, \forall \, \text{bounded } D\subseteq Y , \, \exists \, N=N(\bu,D)\in \bZ_0^+ \, : \, \Phi(n,\bu,D) \subseteq B \,\, \forall n\geq N.
$$
\label{def:uniformly_forward_absorbing}
While, a nonempty compact subset $B \subseteq Y $ is called \emph{uniformly forward absorbing} for a family of inputs $\cV\subseteq \cU $ if
$$
\forall \, \text{bounded } D\subseteq Y , \, \exists \, N=N(D)\in \bZ_0^+ \, :  \, \forall \, \bu \in \cV, \, \Phi(k,\bu,D) \subseteq B \,\, \forall k\geq N.
$$
\end{defn}

Recall Definition \ref{def:posit_invar} of $\cU$-positively invariant set.
Below we prove, under very few assumptions, that for a generic RNN with leaky-integrator neurons \eqref{eq:RNN-map}, if the image of the activation function $\phi$ is $(-L,L)$ then the hypercube $[-L,L]^{N_r}$ of phase space is a uniformly forward absorbing, positively invariant set for all input sequences assuming values in a given compact space.
\begin{prop}
\label{prop:absorbing-set-pullback}
Let us consider a compact subspace $U \subset \bR^{N_i}$ as the set of admissible input values and $\cU := U^{\bZ} $ as the set of admissible input sequences.
Let us consider $ Y:=\bR^{N_r}$ equipped with the Euclidean distance.
For all $\bu = \{ u[k] \}_{k \in \bZ} \in \cU$, the input-driven dynamics of a leaky RNN with feedback of the output \eqref{eq:leaky_rnn}-\eqref{eq:output} are ruled by
\begin{equation}
    \label{eq:cocycle-ESN}
    x[k]  = (1-\alpha)x[k-1] + \alpha \phi( W_r x[k-1] + W_{in} u[k] + W_{fb} \psi(x[k-1]) ).
\end{equation}
If $\phi, \psi$ are upper semi-continuous functions and $\phi$ is non-decreasing with image $(-L,L)$, then $L\cdot I^{N_r} := [-L, L]^{N_r}$ is a $\cU$-positively invariant (see Definition \ref{def:posit_invar}) uniformly forward absorbing set for inputs $\cU$ (see Definition \ref{def:uniformly_forward_absorbing}).
\end{prop}

\proof
The proof is divided in two parts. First we prove that $L\cdot I^{N_r} $ is a $\cU$-positively invariant set, then we show that $L\cdot I^{N_r}$ is a uniformly forward absorbing set.

(i) \emph{Positively invariant.}\\
Let be given an initial condition $  x[0]  $ such that $\|  x[0] \|_{\infty} \leq L$. Thanks to the triangle inequality applied on \eqref{eq:RNN-map}, we have
$$
    \lVert x[1] \rVert_{\infty}    \leq  (1 - \alpha) \| x[0] \|_{\infty} + \alpha \| \phi( W_{r} x[0] + W_{in}u[1] + W_{fb} \psi(x[0]) ) \|_{\infty} \leq  (1 - \alpha) L +   \alpha L = L .\\
$$
Analogously, if at any time step it holds that $ \lVert x[N] \rVert_{\infty} \leq L $, then it will be $ \lVert x[k] \rVert_{\infty} \leq L  , \,\, \forall k \geq N $.

(ii) \emph{Uniformly forward absorbing.}\\
First of all, note that in this framework the universe of possible past input sequences coincides with the universe of possible future sequences. 
In other words, for all $\bu \in \cU$ there exists a $\bv \in \cU$ such that $ v[k] = u[-k-1]  $ for all $k \geq 1$, thus driving the system in pullback sense with the past sequence $\bu^-$ is equivalent to drive the system in forward sense with the future sequence $\bv^+$.
Therefore, a pullback absorbing set of Definition \ref{def:pullback_absorbing} for the family $\cU$ is such if and only if it is a forward absorbing set for the family $\cU$. 
We will prove that
$$
\forall \alpha \in (0,1) , \, \forall x[0]  \in Y , \, \exists N=N(\alpha, x[0]) \, : \, \forall \bu^+ \in \cU^+, \quad \|  x[k] \|_{\infty} \leq L , \,\, \forall k \geq N,
$$
which implies that $L\cdot I^{N_r}$ is a uniformly forward absorbing set for the family $\cU$.

The case of $\alpha = 1 $ brings trivially to the thesis. Thus let us suppose that $\alpha \in (0,1)$.
Let be given the initial condition $ x[0] $, where we assume $\| x[0] \|_{\infty} > L$, otherwise the argument of (i) brings to the thesis.
Note that for all $x[k] $ such that $\| x[k] \|_{\infty} > L$ it holds that 
$$
    \| x[k+1] \|_{\infty} \leq (1- \alpha) \| x[k] \|_{\infty} + \alpha L < \| x[k] \|_{\infty} 
$$
for all $k\geq 0 $.
Now, since $\phi, \psi$ are upper semi-continuous functions and $\phi$ is non-decreasing then the function $ \nu:U\times [-R,R]^{N_r} \longrightarrow \bR_{\geq 0} $ defined as $ \nu(u,x):= \| \phi( W_{r} x + W_{in}u + W_{fb} \psi ( x ) ) \|_{\infty} $ is upper semi-continuous.
Hence, defined $R:=\|x[0]\|_{\infty} $ and since $U$ is compact, there exists a maximum value 
$$
\eta := \max_{u \in U, \,\,  x :  \|x\|_{\infty} \leq R}\| \phi( W_{r} x + W_{in}u + W_{fb} \psi ( x ) ) \|_{\infty}.
$$
Exploiting recursively the triangle inequality on \eqref{eq:cocycle-ESN} the following holds
\begin{align*}
\label{eq:}
    \|  x[k] \|_{\infty}   \quad  &\leq \quad (1 - \alpha) \| x[k-1] \|_{\infty} + \alpha \eta \quad \leq \quad (1 - \alpha)^2 \| x[k-2] \|_{\infty} + \alpha \eta (1-\alpha) + \alpha \eta \quad \leq \\
                                &\leq \quad (1 - \alpha)^k \| x[0] \|_{\infty} + \eta \alpha \sum_{j=0}^{k-1}(1-\alpha)^j \quad =  \quad (1 - \alpha)^k \| x[0] \|_{\infty} + \eta \left[  1 - \alpha \sum_{j=k}^{\infty}(1-\alpha)^j  \right],
\end{align*}
where the last equality holds true in virtue of the geometric series limit $  \alpha \sum_{j=0}^{\infty}(1-\alpha)^j = 1 $.
Now the following inequalities are equivalent,
\begin{align*}
    &(1 - \alpha)^k \| x[0] \|_{\infty} + \eta \Bigl[  1 - \alpha \sum_{j=k}^{\infty}(1-\alpha)^j \Bigl] \leq L  \Longleftrightarrow  (1 - \alpha)^k \| x[0] \|_{\infty} \leq L - \eta     + \eta \alpha \sum_{j=k}^{\infty}(1-\alpha)^j \Longleftrightarrow \\
    & \| x[0] \|_{\infty} \leq \eta \alpha \sum_{j=k}^{\infty}(1-\alpha)^{j-k}  + \dfrac{L -\eta}{(1 - \alpha)^k} \Longleftrightarrow  \| x[0] \|_{\infty} \leq \eta \underbrace{ \alpha \sum_{j=0}^{\infty}(1-\alpha)^{j} }_{=1} + \dfrac{L -\eta}{(1 - \alpha)^k} \Longleftrightarrow\\
    & \| x[0] \|_{\infty} -  \eta \leq  \dfrac{L -\eta}{(1 - \alpha)^k}.
\end{align*}
Note that, by the boundedness hypothesis of $\phi$, it holds that $\eta \leq L $. 
If $\eta = L $ then we conclude from the last inequality that 
$$
    \| x[0] \|_{\infty} \leq  \eta = L,
$$
which is in contradiction with the assumption of $\| x[0] \|_{\infty} > L$.
Accordingly, $\eta < L $ must hold, leading to the following inequality:
\begin{align*}
    & (1 - \alpha)^k \leq \dfrac{L -\eta}{\| x[0] \|_{\infty} -  \eta} \quad \Longleftrightarrow \quad k \geq  \dfrac{ \ln(L - \eta) - \ln( \| x[0] \|_{\infty} -  \eta ) }{\ln(1 - \alpha)}.
\end{align*}
Therefore, after a number of time steps given by $N(\alpha, x[0]):=\dfrac{ \ln(L - \eta) - \ln( \| x[0] \|_{\infty} -  \eta ) }{\ln(1 - \alpha)}$, the internal state of a leaky ESN \eqref{eq:RNN-map} will surely lie inside the hypercube $L\cdot I^{N_r}$.
\qed

For the sake of clarity, we report here below without proof \cite[Theorem 3.20]{kloeden2011nonautonomous} but using our notation. 
\begin{thm}
\label{thm:existence-pullback}\cite[Theorem 3.20]{kloeden2011nonautonomous}
Let $U \subset \bR^{N_i}$ be compact and $(\cU, d_{\cU})$ be the compact metric space of admissible input sequences, where $ \cU:= U^{\bZ} $ and $d_{\cU}$ as defined in \eqref{eq:metric-inputsequences}. Let $\sigma: \cU \longrightarrow \cU$ be the shift operator. 
Let $(\sigma, \Phi)$ be the skew product flow on a complete metric space $(Y,d_Y)$, and $\Phi$ defined as Definition \ref{def:cocycle_mapping}.
If there exists a nonempty compact subset $ B \subset Y $ which is \emph{pullback absorbing} and \emph{positively invariant} for $\cU$, then there exists a unique pullback attractor $\bA= \{ A_{\bu} \}_{\bu \in \cU}$ with fibres in $B$ uniquely determined by
\begin{equation}
    \label{eq:exist-pullback}
    A_{\bu} := \bigcap_{m\geq 0} \overline{ \bigcup_{s \geq m} \Phi(s,\sigma^{-s}(\bu), B) }, \qquad \forall \, \bu \in \cU.
\end{equation}
In addition, since $(\cU, d_{\cU})$ is a compact metric space the subset $A(\cU):= \overline{ \bigcup_{\bu \in \cU}A_{\bu} } \subset B $ uniformly (in $\cU$) attracts every bounded set of the phase space, that is
\begin{equation}
    \label{eq:uniform-limit-attractor}
    \lim_{k \rightarrow \infty} \sup_{\bu \in \cU} h( \Phi(k, \bu, D) , A(\cU) ) = 0, \qquad \forall \text{ bounded } D \subseteq X.
\end{equation}
\end{thm}

Thanks to Proposition \ref{prop:absorbing-set-pullback} we can apply Theorem \ref{thm:existence-pullback} on a RNN taking values with phase space the whole $\bR^{N_r}$ and use the set $ B = [-L,L ]^{N_r} $ in order to construct the pullback attractor.
Moreover, since $\cU$ is compact the second part of Theorem \ref{thm:existence-pullback} implies that the entire nonautonomous dynamics of a RNN is uniformly attracted to a closed subset inside $ B = [-L,L ]^{N_r} $.
This justifies our assumption of considering the whole space as $ X = [-L,L ]^{N_r}$ whenever referring to the RNN nonautonomous dynamics even with leaky neurons.

Therefore, pullback attractor has component sets made by 
\begin{equation}
    \label{eq:pullback-sets}
    A_{\bu} := \bigcap_{m\geq 0} \overline{ \bigcup_{s \geq m} \Phi(s,\sigma^{-s}(\bu), X) } \qquad \forall \, \bu \in \cU.
\end{equation}
Note that under this formalism the resulting set $A_{\bu}$ actually depends only by the left-infinite sequence $\bu^-= \{ \ldots, u[-2], u[-1], u[0] \}$.
Furthermore, since $\cU=U^\bZ$ is shift-invariant, i.e. $\sigma^n(\cU) = \cU$ for all $n \in \bZ$, we can equivalently write \eqref{eq:pullback-sets} as follows
\begin{equation}
    \label{eq:shift-pullback-sets}
    A_{\sigma^n(\bu)} := \bigcap_{m\geq -n} \overline{ \bigcup_{s \geq m} \Phi(s+n,\sigma^{-s}(\bu), X) } =  
    \bigcap_{m'\leq n} \overline{ \bigcup_{s' \leq m'} \Phi(n-s',\sigma^{s'}(\bu), X) }  \qquad \forall \, \bu \in \cU,
\end{equation}
where the last equality is obtained transforming indices as $m'=-m$ and $s'=-s$.

As a consequence, for any fixed input sequence $\bv \in \cU$, from \eqref{eq:shift-pullback-sets} we get the component sets of the pullback attractor with regard to such input sequence:
\begin{equation}
    \label{eq:pullback-sequence}
    \bA(\bv) = \{    A_{\sigma^n(\bv)} \}_{n \in \bZ}  .
\end{equation}
Proposition \ref{prop:pullback_thm} below implies that \eqref{eq:pullback-sequence} is exactly the natural association \eqref{eq:natural_association} defined in Section \ref{sec:cocycle_formalism}, for a given input sequence $\bv \in \cU$.
\begin{lem}
\label{lem:nested-cocycle}
Assume the hypotheses of Theorem \ref{thm:existence-pullback} hold. If there exists a $\cU$-positively invariant set $ B \subset X$, then 
\begin{equation}
    \label{eq:nested-cocycle}
    \Phi(N+1, \bv, B) \subseteq \Phi(N, \sigma(\bv), B) \qquad \forall N \in \bZ_0^+, \, \forall \bv \in \cU.
\end{equation} 
\end{lem}
\proof
The base case of $N=0$, which reads $ \Phi(1, \bv, B) \subseteq B $, holds thanks to the fact that $B$ is $\cU$-positively invariant.
Therefore, the result follows from the \emph{cocycle property} \eqref{eq:cocycle_property},
$$
    \Phi(N+1, \bv, B) = \Phi \bigl( N, \sigma(\bv),  \Phi(1, \bv, B) \bigl) \, \subseteq \, \Phi \bigl( N, \sigma(\bv),  B \bigl).
$$
\qed

\begin{prop}
\label{prop:pullback_thm}
Let $U \subset \bR^{N_i}$ be compact and $(\cU, d_{\cU})$ be the compact metric space of admissible input sequences, where $ \cU:= U^{\bZ} $ and $d_{\cU}$ as defined in \eqref{eq:metric-inputsequences}. 
Let $\sigma: \cU \longrightarrow \cU$ be the shift operator. 
Let $(\sigma, \Phi)$ be the skew product flow on a complete metric space $(X,d_X)$, $ X \subseteq \bR^{N_r}$ closed and $d_X$ the Euclidean distance, and $\Phi$ defined as Definition \ref{def:cocycle_mapping}.
Let us be given an input sequence $\bv \in \cU$, yielding the subset $\cV \subseteq \cU $ of input sequences $ \cV := \{ \ldots, \sigma^{-2}(\bv), \sigma^{-1}( \bv), \bv ,\sigma^{1}( \bv), \sigma^{2}( \bv), \ldots \}  $.
If there exists a nonempty compact subset $ B \subset X$ which is \emph{pullback absorbing} and \emph{positively invariant} for $\cV$, then there exists the (global) pullback attractor $\bA(\bv):= \{ A_n \}_{n \in \bZ}$ of the dynamics driven by $\bv$ and it has component sets 
\begin{equation}
    \label{eq:equivalent-pullback}
     A_n := \bigcap_{m \leq n} \Phi( n-m, \sigma^m(\bv), B ), \qquad  \forall n \in \bZ.
\end{equation}
In addition, if $\cV$ is compact in $(\cU, d_{\cU})$ then $ A(\cV):= \overline{ \bigcup_{n \in \bZ }A_n } \subset B $ uniformly (in time) attracts the whole phase space driven by $\bv$, that is
\begin{equation}
    \label{eq:equivalent-uniform-limit-attractor}
    \lim_{k \rightarrow \infty} \sup_{n \in \bZ} h( \Phi(k, \sigma^n(\bv), X) , A(\cV) ) = 0.
\end{equation}
\end{prop}
\proof
The proof is a direct application of Theorem \ref{thm:existence-pullback}. 
We need to prove that component sets of \eqref{eq:exist-pullback} coincide with component sets of \eqref{eq:equivalent-pullback}.
Let be given a $\bv \in \cU$. Fixed a $n \in \bZ$, let us define a sequence of sets $B_s:= \Phi(n+s,\sigma^{s}(\bv), B)$, for $s\geq -n$. $B_s=\overline{B_s}$ holds as $B$ is compact and $\forall N \in \bZ_0^+, \, \forall \bv \in \cU$ the function $ \Phi(N, \bv, \cdot) ~ : ~ X \longrightarrow ~ X $ is a \emph{closed map}, i.e. it maps closed sets in closed sets\footnote{This is ensured by the \emph{closed map lemma}, which holds since $X$ is a compact space (as domain) and a Hausdorff space (as codomain) and $ \Phi(n, \bu, \cdot) : X \longrightarrow X $ is continuous $\forall n \in \bZ, \, \forall \bu \in \cU$. }, hence compact sets in compact sets in our case.
It is known that, for a nonincreasing sequence of closed sets $\{ B_s \}_{s\geq -n}$, the limit set $ L(n) := \bigcap_{m\geq -n}  \bigcup_{s \geq m} B_s $ exists and it also holds that $B_m = \bigcup_{s \geq m} B_s$. 
Therefore, it is sufficient to prove that $B_{s+1} \subseteq B_s $ for having that $ L(n) := \bigcap_{m\geq -n} B_m$. 
Relation $B_{s+1} \subseteq B_s $ holds thanks to Lemma \ref{lem:nested-cocycle}.
Concluding, component sets \eqref{eq:exist-pullback} can be written as $ L(n) := \bigcap_{m\geq -n} B_m $, which reads 
$$
   A_n := \bigcap_{m\geq -n} \Phi(n+m,\sigma^{m}(\bv), B)=\bigcap_{m'\leq n} \Phi(n-m',\sigma^{m'}(\bv), B) ,
$$
that is \eqref{eq:equivalent-pullback} of thesis. \\
To conclude, \eqref{eq:equivalent-uniform-limit-attractor} follows from \eqref{eq:uniform-limit-attractor} by noting that we are interested in the subset of input sequences given by $ \cV = \{ \ldots, \sigma^{-2}(\bv), \sigma^{-1}( \bv), \bv ,\sigma^{1}( \bv), \sigma^{2}( \bv), \ldots \}  $. 
Therefore, the supremum $ \sup_{ \bu \in \cV} $ can be expressed with the supremum $ \sup_{ n \in \bZ } $ of \eqref{eq:equivalent-uniform-limit-attractor} of thesis.
\qed
\begin{remark}
Note that, in general, for a given $\bv \in \cU $ the set $\cV = \{ \sigma^n(\bv) \}_{ n \in \bZ}$ is not compact in $\cU$.
 In particular, if $\bv$ is aperiodic then $\cV$ will have limit points that are not contained in $\cV$. 
\end{remark}

\subsection{ESP in RNNs implies uniformly (in time) forward convergence} 
\label{sec:forward}

For the sake of completeness, in Theorem \ref{thm:our_uniform_attractor} we show, using our framework, a result found in the literature, which links pullback attraction with forward attraction in the particular case where the global pullback attractor is an entire solution.
\begin{thm}
\label{thm:our_uniform_attractor}\cite[Theorem 3.44]{kloeden2011nonautonomous}  
Let $U \subset \bR^{N_i}$ be compact and $(\cU, d_{\cU})$ be the compact metric space of admissible input sequences, where $ \cU:= U^{\bZ} $ and $d_{\cU}$ as defined in \eqref{eq:metric-inputsequences}. Let $\sigma: \cU \longrightarrow \cU$ be the shift operator. 
Let $(\sigma, \Phi)$ be the skew product flow on the complete metric space $ Y = \bR^{N_r}$, with the Euclidean distance, and $\Phi$ defined as Definition \ref{def:cocycle_mapping}.
Assume there exists a nonempty \emph{compact} subset $ B \subset Y $ such that:
\begin{itemize}
    \item $B$ is a \emph{$\cU$-positively invariant set}; and
    \item $B$ is \emph{uniformly forward absorbing for $\cU$}.
\end{itemize}
Suppose that for all $\bu \in \cU$ the (global) pullback attractor for input $\bu$ is an entire solution.
Then, for any given input sequence $\bv \in \cU$ such entire solution for input $\bv$, denoted as $\bx=\{ x[k] \}_{k \in \bZ}$, is also attracting in forward sense uniformly in time:
$$
    \lim_{k \rightarrow \infty } \sup_{n \in \bZ} h( \Phi( k, \sigma^n(\bv), D ) ,  x[n+k]  ) =0,  \qquad \forall \text{ bounded } D\subseteq Y.
$$
\end{thm}

In Proposition \ref{prop:absorbing-set-pullback} we proved that $B=[-L,L]^{N_r}$ fulfils the hypothesis of Theorem \ref{thm:our_uniform_attractor} for a generic leaky RNNs with $\phi, \psi $ upper semi-continuous functions and $\phi$ non-decreasing with bounded image.
Therefore, Theorem \ref{thm:our_uniform_attractor} for a generic RNN reads as: if the global pullback attractor is an entire solution for all $\bu \in U^\bZ$ then such unique entire solution of the system is also forward attracting uniformly in time.
More formally, we can state the following result.
\begin{prop}
\label{prop:ESP_unif-forw-attract}
Let be given an $\alpha \in (0,1]$ and real matrices $ W_r, W_{in}, W_{fb} $ of dimensions, respectively, $  N_r\times N_r,  \, N_r \times N_i , \,  N_r\times N_o $.
Let $\phi:\bR \rightarrow (-L,L), \,\, \psi: \bR^{N_r} \rightarrow \bR^{N_o} $ be upper semi-continuous functions and $\phi$ non-decreasing.
Consider the following input-driven leaky RNN with feedback of the output
\begin{align}
    \label{eq:inp-driv_RNN}
    & x[k] = G(u[k],x[k-1]), \qquad x \in \bR^{N_r}, \,\, u \in U \subset \bR^{N_i} \text{ compact},\\
    & G(u,x) = (1-\alpha) x + \alpha \phi( W_r x + W_{in} u + W_{fb} \psi(x) ).
\end{align}
If the ESP as originally introduced in \cite[Definition 1]{jaeger2001echo} holds for the input-driven RNN \eqref{eq:inp-driv_RNN} w.r.t the compact input space $\cU=U^\bZ$ then the unique entire solution of such  nonautonomous system is \emph{uniformly state contracting} \cite[Definition 4]{jaeger2001echo}.
\end{prop}

\noindent\textbf{Proof of Lemma \ref{lem:max_sing_value}}\\
\label{proof:lemma}
Let us proceed by contradiction assuming that $\sigma(A) > M $.
By definition,
$$
\sigma(A) = \max_{y \in \bR^{N_r} \setminus{\{0\}}} \dfrac{\|A y \|}{\|y \|}.
$$
Therefore, for a small enough $\varepsilon>0$, it must exist a unit vector $\overrightarrow{v}$ such that
\begin{equation}
    \label{eq:linear_nice_direction}
    \|A (c \overrightarrow{v}) \| \geq (M+\varepsilon) \| c \overrightarrow{v} \| \qquad \forall c \in \mathbb{R} \setminus \{0\}.
\end{equation}
Now let us move on the line $z=x^* + c \overrightarrow{v} $ and consider the linearisation of the map $F$ around $ x^* $, which reads
\begin{equation}
    \label{eq:linearisation}
    F(z) = x^* + A (z-x^* ) + R(z-x^*),
\end{equation}
where the rest of the expansion $R(z - x^*)$ is such that
\begin{equation}
    \label{eq:rest_linearisation}
    \lim_{z \rightarrow{x^*}} \dfrac{ \| R(z-x^*) \| }{ \| z-x^* \| }=0,
\end{equation}
whenever the map $F$ is regular enough in $x^*$.
Now, since $x^*$ is a UASP, there exists a $\delta > 0 $ such that
\begin{equation}
    \label{eq:linear_UASP}
    \| F(z) - x^* \| <  M \| z-x^*\|,   \qquad   \forall z \in B_\delta(x^*).
\end{equation}
By means of the expansion \eqref{eq:linearisation} and applying the reverse triangle inequality we get
\begin{equation}
    \label{eq:linear_inverse_triangle}
    \Bigl| \| A (z-x^* ) \| -  \| R(x-x^*) \| \Bigl| \leq \| A (x-x^*) + R(x-x^*) \| = \| F(x) - x^* \| <  M \| x-x^* \| .
\end{equation}
Now, if $z$ is close enough to $x^*$ along the direction pointed by $ \overrightarrow{v} $, then it holds that $ \| A (z-x^*) \| \geq  \| R(z-x^*) \| $. 
Indeed, \eqref{eq:rest_linearisation} implies that for any $ \epsilon> 0 $ there exists a $c_\epsilon>0$ such that $  \| R(c \overrightarrow{v} ) \|  \leq \epsilon \| c \overrightarrow{v} \| $ for all $|c|\leq c_\epsilon$.
Moreover, exploiting the fact that $z$ is chosen such that $ z- x^* = c\overrightarrow{v} $ satisfies \eqref{eq:linear_nice_direction}, then it holds that $ (M+\varepsilon) \| c \overrightarrow{v}\|  =  (M+\varepsilon)  \| z- x^*\| \leq \| A (z-x^*)\| $.
Hence, for the choice $ \epsilon= (M+\varepsilon) $ there exists a $c_{M+\varepsilon} > 0 $ such that $ \| R( z - x^* ) \| \leq \| A ( z - x^* )\| $ for all $ z - x^* = c \overrightarrow{v} $ having $|c|\leq c_{M+\varepsilon} $.
Therefore, we can rewrite \eqref{eq:linear_inverse_triangle} as follows
\begin{equation}
    \label{eq:linear_inverse_triangle_2}
    \| R(z-x^*) \|  > \| A (z-x^* ) \|   -  M \| z-x^* \| .
\end{equation}
Finally, applying again \eqref{eq:linear_nice_direction} in \eqref{eq:linear_inverse_triangle_2}, we obtain 
\begin{equation}
    \label{eq:linear_last_inequality}
    \| R(z-x^*) \|  >  ( M+\varepsilon  -  M ) \| z-x^* \| ,
\end{equation}
i.e. for all $ z - x^* = c \overrightarrow{v} $ having $|c|\leq c_{M+\varepsilon} $ holds $  \| R(z-x^*) \|  > \varepsilon \| z-x^* \| $, which is in contradiction with relation \eqref{eq:rest_linearisation}. 
\qed

\section{Training details for the context-dependent task in Section~\ref{sec:Maass}}
\label{sec:training}

We trained an RNN of the form \eqref{eq:leaky_rnn}-\eqref{eq:output} with $N_r=500$ neurons by means of a supervised learning algorithm (ridge regression) which exploits the output feedback as a mechanism for optimising the recurrent layer; see \cite{ceni-ff-18} for more details.
The state-update \eqref{eq:leaky_rnn} was configured with $\alpha = 1$.
The readout function in \eqref{eq:output} was parametrised through a matrix $W_o$ so that $\psi(x) = W_o  x$. 
In the training phase, Gaussian noise was added inside the activation function $\phi$, see \cite[Equation 5]{ceni-ff-18}, with zero mean and standard deviation set to $0.05$. 
The entries of matrices $W_{in} , W_{fb}$ and $W_r$ were i.i.d. drawn from a uniform distribution in $[-1, 1]$; the sparseness of $W_r$ was set to $95\%$. 
Moreover, the matrix $W_r$ was rescaled so that its spectral radius equals to $0.9$. 
Finally, the readout matrix $W_o$ have been determined via ridge regression, with regularisation parameter set to $ 0.7 $.
We consider input sequences with $15000$ time-steps.
The first $10000$ time-steps have been used for training and the remaining $5000$ for testing.
Once the training session was completed, we closed the feedback loop by injecting the output $z_1[k]$ into the network state-update.

\end{document}